\newif\ifarxiv
\arxivtrue
\pdfoutput=1

\newcommand{\de}{{d}}
\newcommand{\e}{\operatorname{e}}
\newcommand{\p}{\partial}
\newcommand{\ney}{\boldsymbol{r}'}
\newcommand{\nex}{\boldsymbol{r}}

\newcommand{\nephi}{\boldsymbol{\phi}}

\newcommand{\nn}{\boldsymbol{n}}

\newcommand{\densi}{\eta}
\newcommand{\inc}{\text{inc}}
\newcommand{\scat}{\text{scat}}
\newcommand{\gammainc}{\Gamma_j^\inc}
\newcommand{\gammainctilde}{\widetilde{\Gamma}_j^\inc}
\newcommand{\gammainf}{{\Gamma}_j^\infty}

\ifarxiv

  \documentclass[11pt,letter]{article}

  \usepackage[hidelinks=true]{hyperref}
  \usepackage{amsmath}
  \usepackage{amsfonts}
  \usepackage{amssymb}
  \usepackage{amsthm}
  \usepackage{graphicx}
  \usepackage{cite}
  \usepackage{empheq}
  \usepackage{caption}
  \usepackage{enumerate}
  \usepackage{subfig}
  \usepackage{mathrsfs} 

  \newtheorem{remark}{Remark}[section]

\else

\fi

\begin{document}

\title{Windowed Green Function Method for\\ Nonuniform Open-Waveguide
  Problems}

\ifarxiv

\author{Oscar P.
  Bruno$^1\footnote{E-mail:
    \texttt{\{obruno,emmanuel.garza\}@caltech.edu, cperezar@mit.edu}}$,
  Emmanuel Garza$^1$, Carlos P\'erez-Arancibia$^{2}$\\ \\
  \small{$^1$Computing \& Mathematical Sciences, California Institute of
    Technology}\\
  \small{$^2$Department of Mathematics, Massachusetts Institute of Technology}}
       
  \date{\today}

\else

\fi

\maketitle

\begin{abstract}

  This contribution presents a novel Windowed Green Function (WGF)
  method for the solution of problems of wave propagation, scattering
  and radiation for structures which include open (dielectric)
  waveguides, waveguide junctions, as well as launching and/or
  termination sites and other nonuniformities. Based on use of a
  ``slow-rise'' smooth-windowing technique in conjunction with
  free-space Green functions and associated integral representations,
  the proposed approach produces numerical solutions with errors that
  decrease faster than any negative power of the window size. The
  proposed methodology bypasses some of the most significant
  challenges associated with waveguide simulation. In particular the
  WGF approach handles spatially-infinite dielectric waveguide
  structures without recourse to absorbing boundary conditions, it
  facilitates proper treatment of complex geometries, and it
  seamlessly incorporates the open-waveguide character and associated
  radiation conditions inherent in the problem under
  consideration. The overall WGF approach is demonstrated in this
  paper by means of a variety of numerical results for two-dimensional
  open-waveguide termination, launching and junction problems.
\end{abstract}


\section{Introduction} \label{sec:intro}

This paper considers the problem of evaluation of wave propagation and
scattering in nonuniform open-waveguide structures.  This is a problem of
fundamental importance in a wide range of areas, including modeling and design
of dielectric antenna systems, photonic and optical devices, dielectric RF
transmission lines, etc. The numerical simulation of such structures presents
significant challenges---in view of the unbounded character of the associated
dielectric boundaries and propagation domains as well as the presence of
radiating fields, inhomogeneities, and scattering obstacles.

The present contribution introduces an effective methodology for the
solution of such nonuniform open-waveguide problems. Based on use of
Green functions and integral equations akin to those used in the
Method of Moments~\cite{Rao1982}, and incorporating as a main novel
element a certain ``slow-rise'' windowing function, the proposed
Windowed Green Function approach (WGF) can be used to model, with
high-order accuracy, highly-complex waveguide structures---without
recourse to use of mode matching (which can be quite challenging in
the open-waveguide context), absorbing boundary conditions,
staircasing or time-domain simulations.

The finite-difference time-domain method (FDTD) is one of the simplest
and most reliable extant methods for solution of open-waveguide
problems. In the FDTD approach, unbounded domains are truncated by
relying on absorbing boundary conditions or absorbing layers such as
the PML~\cite{Berenger1994}. Further, subpixel smoothing
techniques~\cite{Farjadpour2006} are often used in FDTD
implementations to model material interfaces while maintaining second
order accuracy---in spite of the staircasing that accompanies
Cartesian discretization of curved boundaries. In order to obtain the
frequency response from a FDTD simulation, finally, Fourier transforms
in time are typically used. In spite of its usefulness, the FDTD
approach does present a number of difficulties in the context of
waveguide problems~\cite[p. 223]{taflove} concerning (i)~Illumination
by specified waveguide modes and inadvertent excitation of unwanted
modes; (ii)~Necessary use of sufficiently large computational domains
to allow decay of reactive fields; (iii)~Need for substantially
prolonged simulation times in order for spectral energy above the
waveguide cutoff frequency to reach a given interaction structure of
interest; and (iv)~Necessary use of fine spatial and temporal
discretizations to mitigate the numerical dispersion associated with
the second-order accuracy of the method.

A few Green function methods for open-waveguide problems have also
been proposed. The recent boundary-element method~\cite{Zhang2011},
for example, in which a conductive absorber is used to truncate the
unbounded waveguide structure, requires the excitation source to lie
within the computational domain. The source then produces a radiation
field that decays as $\mathcal{O}(1/r)$ and thus limits the accuracy
of the implementation. Other approaches for open-waveguide problems
include perturbation methods~\cite{Ciraolo2008,Ciraolo2005} which can
effectively handle limited types of (sufficiently small) localized
inhomogeneities.

Relying on the free space Green function and associated integral
equations along the dielectric boundaries, the WGF method presented in
this paper utilizes a slow-rise window function to truncate the
infinite integrals while providing super-algebraically small error
(that is, errors smaller than any negative power of the window
size). The method can easily incorporate bound modes and arbitrary
beams as illuminating sources, and it can treat general
inhomogeneities and multiple arbitrarily oriented waveguides without
difficulty.

The proposed use of slow-rise window functions has been previously
found highly effective in the contexts of scattering by periodic rough
surfaces~\cite{Monro2007,Bruno2014,royal} and obstacles in presence of
layered-media~\cite{Bruno2016,Perez-Arancibia2017} as well as
long-range volumetric propagation~\cite{chaubell2009}. The
implementation details vary from problem to problem; in the present
open-waveguide context, for example, the treatment of incident fields
and windowed integral operators differs significantly from those used
previously.

This paper is organized as follows: Notations and mathematical
background on the open-waveguide problem are presented in
Sec.~\ref{sec:math_form}. The WGF method is then introduced in
Sec.~\ref{sec:wgf}, which includes, in particular, an integration
example which demonstrates in a very simple context the properties of
the slow-rise windowing function. A variety of applications of the
open-waveguide WGF method are presented in Sec.~\ref{numer},
demonstrating applicability to waveguide junctions (including couplers
and sharp bends), dielectric antennas, focused beam illumination, and,
for reference, an unperturbed waveguide for which the exact solution
is known. In all cases the WGF method provides high accuracies in
computing times of the order of seconds. Significant additional
acceleration could be incorporated on the basis of equivalent sources
and Fast Fourier Transforms, along the lines of
reference~\cite{royal}; such acceleration methods are not considered
here in view of the fast performance that the unaccelerated method
already provides in the present two-dimensional
case. Sec.~\ref{sec:conclusions}, finally, highlights the properties
of the method, and presents a few concluding remarks.

\ifarxiv

\begin{figure} \centering \subfloat[Schematic depiction of a
  typical waveguide scattering problem.]
  {\includegraphics[width=3.1in]{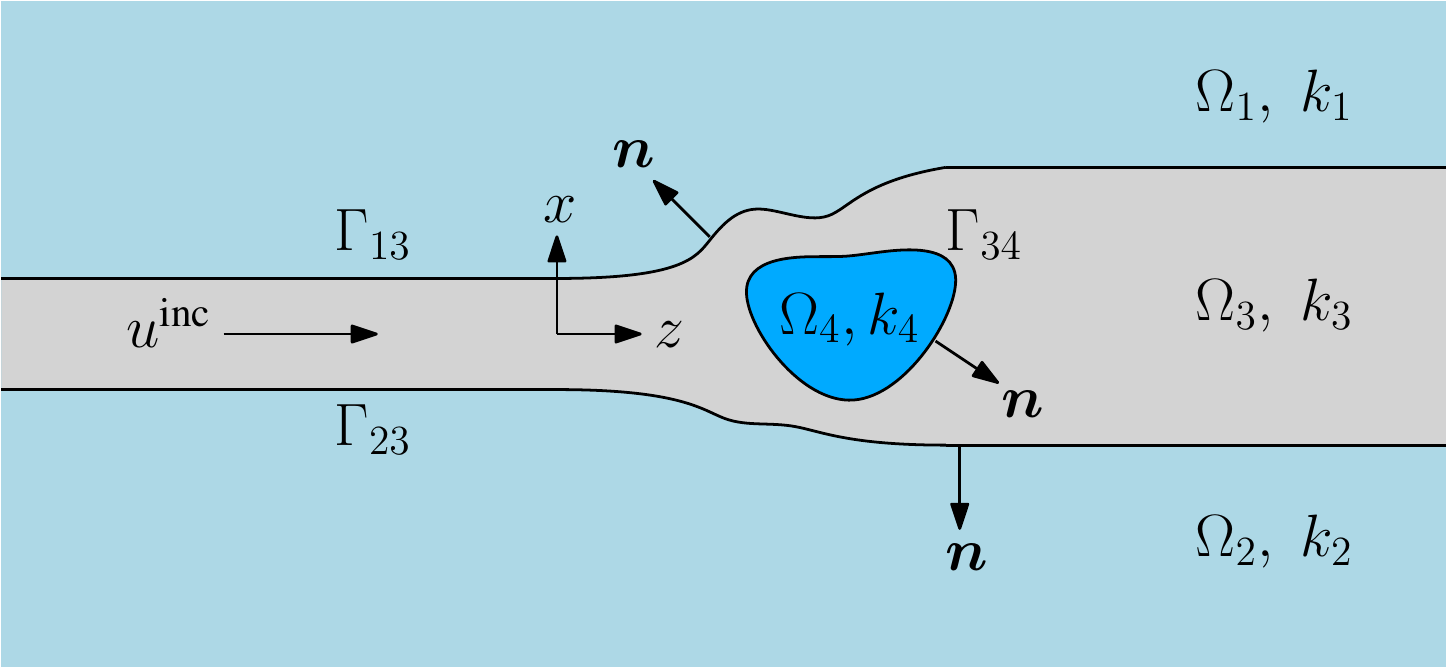}%
    \label{fig:phys}} \quad \subfloat[Notations associated with this
  paper's WGF method.]
  {\includegraphics[width=3.1in]{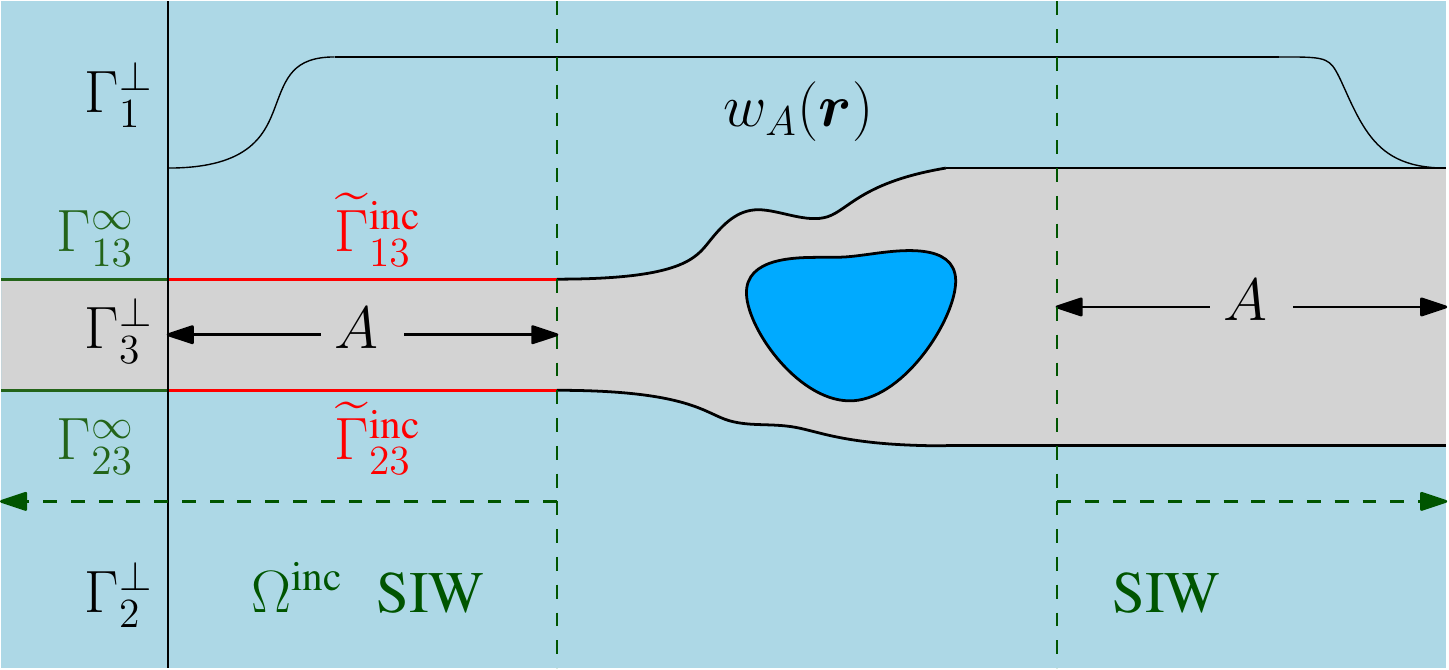}%
    \label{fig:wgf}} 
  \caption{The open-waveguide problem and geometrical structures
    utilized in the WGF method presented in this paper.}
  \label{fig:domains}
\end{figure}

\else


\fi

\section{Mathematical Framework for 2D Open
Waveguides} \label{sec:math_form}

This paper considers the problem of electromagnetic wave propagation and
scattering induced by two-dimensional ($zx$-plane) nonuniform open (dielectric)
waveguides, with application to nonuniformities such as waveguide junctions,
illumination and termination regions. In Fig.~\ref{fig:domains} a schematic
depiction of the problem is presented.

General two-dimensional structures consisting of spatial arrangements
of two-dimensional dielectric waveguides and bounded dielectric
structures can be considered within the proposed
framework---including, for example, configurations which are
constructed as a combination of a given finite number of
``semi-infinite waveguide'' structures (SIW) and additional bounded
dielectric bodies. Here a SIW is one of the two portions that result
as a fully uniform waveguide is cut by a straight line (plane)
orthogonal to the waveguide axis. A simple such configuration is
depicted in Fig.~\ref{fig:domains}.

Additionally, it is useful to identify {\em connected} (bounded or unbounded)
regions that are occupied by a given dielectric material; as illustrated in
Fig.~\ref{fig:domains}, in this paper these regions are denoted by $\Omega_j$
($j=1,\dots,N$). The electrical permittivity, magnetic permeability, refractive
indices and wavenumbers in $\Omega_j$ are denoted $\varepsilon_j$, $\mu_j$,
$n_j$ and $k_j=\omega \sqrt{\varepsilon_j\mu_j}=n_j\omega/c$ ($c$ being the
speed of light in vacuum), respectively.

The structure may be illuminated by arbitrary combinations of bound
waveguide modes supported on a single component SIW (or, more
precisely, by the restriction to the given SIW of a mode of the
corresponding fully uniform waveguide). By linearity, simultaneous
illumination by several SIWs can be obtained directly by addition of
the corresponding solutions for single SIW illumination. Letting $\nex
= (z,x)$ and denoting by $\mathbf{\chi}^{\inc}$ the indicator function
of a SIW region $\Omega^\inc$ that contains the prescribed
illuminating field,
\begin{equation} \label{indicator} 
\mathbf{\chi}^{\inc} (\nex) =
\begin{cases} 1 & \text{for }\nex \in \Omega^\inc \\ 0 &
\text{for }\nex \not\in \Omega^\inc,
\end{cases}
\end{equation} 
the total electric field $\boldsymbol{E}$ is given by $\boldsymbol{E}
= \boldsymbol{E}^\inc \mathbf{\chi}^\inc+\boldsymbol{E}^\scat$.  (See
Fig.~\ref{fig:wgf}, where the region $\Omega^\inc = \{ z<0 \}$ is such
that the waveguide boundaries are flat within $\Omega^\inc$, and thus,
a bound mode can be prescribed in this region as an incident field.)
The scattered field $\boldsymbol{E}^\scat$ is assumed to satisfy an
appropriate radiation condition (which, roughly, states that the
scattered field propagates away from all inhomogeneities as either
outward waveguide modes or cylindrical waves;
see~\cite[Eq. (24)]{Nosich1994} for details) in each component
$\Omega_j$ that is not bounded ---in addition to the Maxwell equations
which, in the two-dimensional case considered in what follows, reduce
to the Helmholtz equation.

Assuming a time-harmonic temporal dependence of the form $\e^{-i\omega
  t}$ (which is suppressed in all expressions in this paper) and
letting $u$ (resp.  $u^\scat$) denote either the $y$-component of the
total (resp. scattered) electric field in TE-polarization or the
$y$-component of the total (resp. scattered) magnetic field in
TM-polarization, the field component $u = u^\inc\chi^\inc + u^\scat$
is the unique~\cite{Nosich1994} radiating solution of the problem
\begin{equation}
  \label{eq:helmholtz1} \left\{ \begin{array}{rcccc} \Delta u + k_j^2
      u &=& 0 &\text{in}& \Omega_j,\smallskip \\ u_+ - u_- &=&
      0 & \text{on}& \Gamma_{j\ell} \;\;  (j <\ell), \medskip\\ 
      \displaystyle\frac{\partial
        u_+}{\partial \nn} - {\nu_{j\ell}}\frac{\partial \displaystyle
        u_-}{\partial \nn} &=& 0 & \text{on}& \Gamma_{j\ell} \;\; (j <\ell).
\end{array}\right. 
\end{equation} 
Here $\nu_{j\ell} = 1$ in TE-polarization and $\nu_{j\ell} =
(k_j/k_\ell)^2$ in TM-polarization; for each pair $(j,\ell)$ with $j
<\ell$, $\Gamma_{j\ell}$ denotes the boundary between $\Omega_j$ and
$\Omega_\ell$; for $\nex\in\Gamma_{j\ell}$, $\nn=\nn(\nex)$ denotes
the unit normal vector to $\Gamma_{j\ell}$ which points into the
``plus side'' $\Omega_j$ of $\Gamma_{j\ell}$ (the plus side of
$\Gamma_{j\ell}$ is defined by the aforementioned condition $j
<\ell$); and for $\nex\in\Gamma_{j\ell}$, finally, \ifarxiv
\begin{align}
  u_
\pm (\nex)=\lim_{\delta\to
    0^+}\big[u(\nex\pm\delta
  \nn(\nex))\big]\quad \mbox{and} \quad
  \frac{\p u_\pm}{\p \nn}(\nex) =\lim_{\delta\to 0^+} \big[\nabla u(\nex \pm
\delta \nn(\nex))\cdot \nn(\nex) \big]\label{eq:bv}.
\end{align}

\else
\begin{align}
  u_{\pm}(\nex) &= \lim_{\delta\to 0^+}u(\nex\pm\delta \nn(\nex)), \\
 \frac{\p u_\pm}{\p \nn}(\nex) &= \lim_{\delta\to 0^+} \nabla u(\nex\pm\delta \nn(\nex))\cdot \nn(\nex). \label{eq:bv}
\end{align}
\fi 
\begin{remark}\label{rem_2.1}  A few comments concerning notations are in
  order: (i)~$\Gamma_{j\ell}$ may be empty for a number of pairs
  $(j,\ell)$: for example $\Gamma_{12}=\emptyset$ for the geometry
  displayed in Fig.~\ref{fig:phys}, and $\Gamma_{j\ell}$ is
  necessarily empty, by definition, whenever $j\geq\ell$.  (ii)~The
  use of indicator functions~\eqref{indicator} makes it possible to
  conveniently specify incident fields in an adequately selected SIW
  region $\Omega^\inc$ (which equals the half plane containing the SIW
  that supports the incident field).  (iii)~Once $u$ is determined by
  solving~\eqref{eq:helmholtz1}, the total electromagnetic field in
  the domain $\Omega_j$ ($j=1,2\dots N$) is given by $\boldsymbol{E}
  =(0,u,0)$, $\boldsymbol{H}=\frac{i}{\omega \mu_j}\left(
    \frac{\partial u}{\partial z},0,- \frac{\partial u}{ \partial
      x}\right)$ in TE-polarization, and
  $\boldsymbol{E}=\frac{i\omega\mu_j}{k_j^2}\left(-\frac{\partial
      u}{\partial z},0,\frac{\partial u}{\partial x}\right)$,
  $\boldsymbol{H}=\left(0,u,0\right)$ in TM-polarization.
\end{remark}

\section{Windowed Green Function Method (WGF)} \label{sec:wgf}


\subsection{Integral Equation Formulation} 
\label{sec:ieqf}
This section presents an integral equation formulation
for the propagation and scattering problem considered in
Sec.~\ref{sec:math_form}. For simplicity, it is assumed that the
structure is illuminated by means of an arbitrary superposition of
bound modes incoming from a single SIW whose optical axis coincides
with the $z-$axis; the generalization to structures containing
multiple arbitrarily-oriented waveguides is straightforward. Under
this assumption, the incident field is prescribed by
\begin{align} \label{eq:inc_mode} u^\inc(z,x) = \sum_{m=1}^M A^\inc_m
  u^m_\perp(x)\e^{ik^m_z z},
\end{align}
where $M$ is the total number of bound modes supported by the
waveguide structure, $A^\inc_m$ denotes the $m$-th modal coefficient,
$u^m_\perp(x)$ is the transverse profile of the mode, and $k^m_z$ is
the corresponding propagation constant for the $m$-th mode. Note that
$u^m_\perp(x)$ and $k^m_z$ can be easily found by solving a one
dimensional eigenvalue problem by means of the method of separation of
variables~\cite{Magnanini2001}. For example, the bound mode solutions
for a single waveguide centered at $x=0$, with half-width $h$ and with
core and cladding wavenumbers $k_\text{co} > k_\text{cl}$
respectively, are given by~\cite{Magnanini2001}
\begin{align} \label{eq:mode}
   u_\perp^m(x) = 
  \begin{cases}
     v\left(\gamma_\text{co}~h \right) e^{-\gamma_\text{cl}(x-h)} &,\quad x > h\\
     v\left(\gamma_\text{co}~x \right) & ,\quad |x|\le h \\
     v\left(-\gamma_\text{co}~h \right) e^{\gamma_\text{cl}(x+h)} &,\quad x < -h,\\
  \end{cases}
\end{align}
where $\gamma_\text{co}=\sqrt{k_\text{co}^2-{(k_z^m)}^2 }$ and
$\gamma_\text{cl}=\sqrt{{(k_z^m)}^2-k_\text{cl}^2}$. Here
$v(\tau)=\cos(\tau)$ for the symmetric modes, $v(\tau)=\sin(\tau)$ for
the antisymmetric modes, and $k^m_z$ ($k_\text{cl}<k^m_z<k_\text{co}$)
is the $m$-th solution of the transcendental equation
\begin{align}
  \begin{cases}
    \displaystyle \gamma_\text{cl} = \nu_\text{wg}
    \gamma_\text{co}\tan{\left(\gamma_\text{co}h\right)}  & \text{(symmetric mode),} \\
    \displaystyle \nu_\text{wg} \gamma_\text{co} = -\gamma_\text{cl}
    \tan{\left(\gamma_\text{co}h\right)} &
    \text{(antisymmetric mode),}
  \end{cases}
\end{align}
with $\nu_\text{wg}=1$ for TE-polarization and
$\nu_\text{wg} = (k_\text{cl}/k_\text{co})^2$ for TM-polarization. As indicated
in Sec.~\ref{numer}, further, additional incident fields such as plane waves,
and finite beams can also be incorporated easily in this context.

In order to introduce the desired system of integral equations, let
$\Gamma_j = \left(\bigcup_{\ell=1}^{j-1}\Gamma_{\ell j}\right) \bigcup
\left(\bigcup_{\ell=j+1}^N\Gamma_{j\ell}\right) $
denote the boundary of the domain $\Omega_j$ and let
$\Gamma=\bigcup_{j=1}^N \Gamma_{j}$ denote the union of all domain
boundaries. Then, calling
\begin{align} \label{eq:densities}
  \varphi(\nex) \equiv u_+(\nex) \quad \text{and} \quad \psi(\nex) \equiv
  \frac{\p u_+}{\p \nn}(\nex), \quad  \nex \in \Gamma,
\end{align}
and using Green's theorem in a manner akin
to~\cite{DeSanto1997} together with the boundary conditions
in~\eqref{eq:helmholtz1}, the representation formula
\begin{equation} \label{eq:rep} 
  u(\nex) = \mathscr{D}[\beta_j \varphi](\nex) - 
  \mathscr{S}\left[\beta_j \nu_j^{-1}\psi\right](\nex),  \quad
  \nex \in \Omega_j
\end{equation}
results, where 
\ifarxiv
\begin{align}
  \beta_j(\nex) =
  \begin{cases}
    \;\;\;1,& \text{for } \nex \in \Gamma_{j \ell} \; (j<\ell) \\
    -1,& \text{for } \nex \in \Gamma_{\ell j} \; (j > \ell)
  \end{cases}, \quad &\quad
  \nu_j(\nex) =
  \begin{cases}
    1,& \text{for } \nex \in \Gamma_{j \ell} \; (j<\ell)  \\
    \nu_{\ell j},& \text{for } \nex \in \Gamma_{\ell j} \; (j>\ell) 
  \end{cases},
\end{align}
\begin{equation}
  \nu(\nex) = \nu_{\ell j}, \quad \text{for } \nex \in \Gamma_{j\ell} \; (j<\ell),
\end{equation}
\else
\begin{align}
  \beta_j(\nex) &=
  \begin{cases}
    \;\;\;1,& \text{for } \nex \in \Gamma_{j \ell} \; (j<\ell) \\
    -1,& \text{for } \nex \in \Gamma_{\ell j} \; (j > \ell)
  \end{cases}, \\
  \nu_j(\nex) &=
  \begin{cases}
    1,& \text{for } \nex \in \Gamma_{j \ell} \; (j<\ell)  \\
    \nu_{\ell j},& \text{for } \nex \in \Gamma_{\ell j} \; (j>\ell) 
  \end{cases},
\end{align}
\begin{equation}
  \nu(\nex) = \nu_{\ell j}, \quad \text{for } \nex \in \Gamma_{j\ell} \; (j<\ell),
\end{equation}
\fi
and where, letting $G_j(\nex,\ney) = \frac{i}{4}H_0^{(1)}(k_j|\nex -
\ney|)$ denote the free-space Green function for the Helmholtz
equation with wavenumber $k_j$, the single and double layer potentials
$\mathscr{S}$ and $\mathscr{D}$ for a given density $\densi$ defined in
$\Gamma$ are given by
\ifarxiv
\begin{equation}\label{eq:layer_potentials}
  \mathscr{S}[\densi](\nex) = 
  \int_{\Gamma_j} G_j(\nex, \ney) \densi(\ney) \de s_{\ney} \quad \text{and} \quad                                                                              
  \displaystyle \mathscr{D}[\densi](\nex) = 
  \int_{\Gamma_{j}} \frac{\partial G_j}{\partial \nn(\ney)}(\nex, \ney) 
  \densi(\ney) \de s_{\ney}, \quad \text{for }\nex \in \Omega_j.
\end{equation}
\else
\begin{align}
  \mathscr{S}[\densi](\nex) &= 
  \int_{\Gamma_j} G_j(\nex, \ney) \densi(\ney) \de s_{\ney},&~\nex \in \Omega_j, \nonumber\\                                                                              
  \displaystyle \mathscr{D}[\densi](\nex) &= 
  \int_{\Gamma_{j}} \frac{\partial G_j}{\partial \nn(\ney)}(\nex, \ney) 
  \densi(\ney) \de s_{\ney},&~\nex \in \Omega_j.
\label{eq:layer_potentials}
\end{align}
\fi
The densities $\varphi$ and $\psi$ in the representation
formula~\eqref{eq:rep}, which, in view of
equation~\eqref{eq:densities}, are given in terms of the total field,
can be expressed as a sum of their incident and scattered
components. In other words $\varphi = \varphi^\inc + \varphi^s$ and
$\psi = \psi^\inc + \psi^\scat$ where, using the indicator
function~\eqref{indicator},
\begin{align} \label{eq:densities2}
  \varphi^\inc(\nex) \equiv \chi^\inc u^\inc_+(\nex)~~&\text{and}~~\psi^\inc(\nex) \equiv
  \chi^\inc \frac{\p u^\inc_+}{\p \nn}(\nex),~~\nex \in \Gamma, \\
  \varphi^\scat(\nex) \equiv u^\scat_+(\nex)~~&\text{and}~~\psi^\scat(\nex) \equiv
  \frac{\p u^\scat_+}{\p \nn}(\nex),~~\nex \in \Gamma.
\end{align}

The desired integral equations for the unknown densities
$\varphi^\scat$ and $\psi^\scat$ are expressed in terms of certain
free-space Green functions and various associated integral operators.
The particular Green function used in the definition of each one of
these operators depends on $\nex$: for $\nex\in\Gamma_{j\ell}$ a
``plus'' (resp. ``minus'') operator uses the Green function $G_j$
(resp. $G_\ell$) corresponding to the refractive index on the plus
side (resp. minus side) of $\Gamma_{j\ell}$. To streamline the
notations in this context, for $\nex\in \Gamma$ let $\Gamma^\pm(\nex)$
be defined as follows: if $\nex\in \Gamma_{j\ell}$ then
$\Gamma^+(\nex) = \Gamma_j$ and $\Gamma^-(\nex)=\Gamma_\ell$
(cf. Remark~\ref{rem_2.1}). The aforementioned integral operators are
thus defined by
\ifarxiv
\begin{equation}\label{eq:opers}
  \begin{array}{r c l r c l}
    \displaystyle S^\pm [\densi](\nex) &=& \displaystyle \int_{\Gamma^\pm(\nex)} G^\pm(\nex,
   \ney) \densi(\ney) ds_{\ney}  , &  
   \displaystyle D^\pm [\densi](\nex) &=& 
   \displaystyle \int_{\Gamma^\pm(\nex)} \frac{\partial G^\pm(\nex, \ney)}{\partial \nn(\ney)} \densi(\ney) ds_{\ney}, \\[.5cm]
   \displaystyle K^\pm [\densi](\nex) &=& 
   \displaystyle \int_{\Gamma^\pm(\nex)} \frac{\partial G^\pm(\nex, \ney)}{\partial \nn(\nex)} \densi(\ney) ds_{\ney} ,  &  
   \displaystyle N^\pm [\densi](\nex) &=& 
   \displaystyle \int_{\Gamma^\pm(\nex)} \frac{\partial^2 G^\pm(\nex, \ney)}{\partial \nn(\nex) \partial \nn(\ney)} \densi(\ney) ds_{\ney}  
  \end{array}  (\nex \in \Gamma).
\end{equation}
\else
\begin{align}\label{eq:opers}
\begin{split}
    \displaystyle S^\pm [\densi](\nex) = \displaystyle \int_{\Gamma^\pm(\nex)} G^\pm(\nex,
   \ney) \densi(\ney) ds_{\ney},~\nex \in \Gamma, \\  
   \displaystyle D^\pm [\densi](\nex) = 
   \displaystyle \int_{\Gamma^\pm(\nex)} \frac{\partial G^\pm(\nex, \ney)}{\partial \nn(\ney)} \densi(\ney) ds_{\ney},~\nex \in \Gamma, \\
   \displaystyle K^\pm [\densi](\nex) =
   \displaystyle \int_{\Gamma^\pm(\nex)} \frac{\partial G^\pm(\nex, \ney)}{\partial \nn(\nex)} \densi(\ney) ds_{\ney} ,~\nex \in \Gamma,  \\  
   \displaystyle N^\pm [\densi](\nex) = 
   \displaystyle \int_{\Gamma^\pm(\nex)} \frac{\partial^2 G^\pm(\nex, \ney)}{\partial \nn(\nex) \partial \nn(\ney)} \densi(\ney) ds_{\ney},~  
  \nex \in \Gamma.
\end{split}
\end{align}
\fi

As is known (cf.~\cite{Colton2013} theorems 3.1 and 3.2), the layer
potentials~\eqref{eq:layer_potentials} satisfy the jump conditions at
$\nex \in \Gamma$:
\ifarxiv
\begin{equation}
  \begin{array}{r c l r c l}
   \displaystyle \lim_{\delta\to 0^+} \mathscr{S}[\densi](\nex\pm\delta \nn (\nex)) &=& S^\pm [\densi](\nex),  &  
   \displaystyle \lim_{\delta\to 0^+} \mathscr{D} [\densi](\nex\pm\delta \nn (\nex)) &=& 
   \displaystyle \pm \frac{1}{2} \densi(\nex) + D^\pm [\densi] (\nex), \\[.5cm]
   \displaystyle \lim_{\delta\to 0^+} \frac{\partial}{\partial \nn}\mathscr{D}[\densi](\nex\pm\delta \nn (\nex)) &=& N^\pm [\densi](\nex),  &
   \displaystyle \lim_{\delta\to 0^+} \frac{\partial}{\partial \nn} \mathscr{S}[\densi](\nex\pm\delta \nn (\nex)) &=& 
   \displaystyle \mp \frac{1}{2} \densi(\nex) + K^\pm [\densi] (\nex).
  \end{array}
\end{equation}
\else
\begin{align}
\begin{split}
   \displaystyle \lim_{\delta\to 0^+} \mathscr{S}[\densi](\nex\pm\delta \nn (\nex)) &= S^\pm [\densi](\nex), \\  
   \displaystyle \lim_{\delta\to 0^+} \mathscr{D} [\densi](\nex\pm\delta \nn (\nex)) &= 
   \displaystyle \pm \frac{1}{2} \densi(\nex) + D^\pm [\densi] (\nex), \\
   \displaystyle \lim_{\delta\to 0^+} \frac{\partial}{\partial \nn}\mathscr{D}[\densi](\nex\pm\delta \nn (\nex)) &= N^\pm [\densi](\nex), \\
   \displaystyle \lim_{\delta\to 0^+} \frac{\partial}{\partial \nn} \mathscr{S}[\densi](\nex\pm\delta \nn (\nex)) &= 
   \displaystyle \mp \frac{1}{2} \densi(\nex) + K^\pm [\densi] (\nex).
\end{split}
\end{align}
\fi
Thus, adding the limits of the fields~\eqref{eq:rep} (resp. the normal
derivatives of the fields~\eqref{eq:rep}) on the plus and minus sides
of $\Gamma$, the system of integral equations
\ifarxiv
\begin{equation} \label{eq:system1} E(\nex) \nephi^\scat (\nex) + T[\nephi^\scat] (\nex)
  = - E(\nex) \nephi^\inc (\nex) - T[\nephi^\inc] (\nex), \quad
  \text{for } \nex \in \Gamma,
\end{equation}
\else
for $\nex \in \Gamma$
\begin{equation} \label{eq:system1} E(\nex) \nephi^\scat (\nex) + T[\nephi^\scat] (\nex)
  = - E(\nex) \nephi^\inc (\nex) - T[\nephi^\inc] (\nex),
\end{equation}
\fi
results, where
\ifarxiv
\begin{equation} \label{eq:oper}
  \displaystyle E(\nex) = \text{diag}\left[1, \frac{1+\nu(\nex)}{2\nu(\nex)}\right], \quad
  T = 
  \begin{bmatrix}
    \displaystyle D^- - D^+ & \displaystyle S^+ - (1/\nu) S^- \\
    \displaystyle N^- - N^+ & \displaystyle K^+ - (1/\nu) K^-
  \end{bmatrix},  
\end{equation}
\else
$\displaystyle E(\nex) = \text{diag}\left[1, \frac{1+\nu(\nex)}{2\nu(\nex)}\right]$ and
\begin{equation} \label{eq:oper}  
  T = 
  \begin{bmatrix}
    \displaystyle D^- - D^+ & \displaystyle S^+ - (1/\nu) S^- \\
    \displaystyle N^- - N^+ & \displaystyle K^+ - (1/\nu) K^-
  \end{bmatrix},  
\end{equation}
\fi
and where the density vectors are given by
\begin{align}
 \nephi^\scat = [\varphi^\scat, \psi^\scat]^t, \quad \text{and} \quad
  \nephi^\inc = [\varphi^\inc, \psi^\inc]^t. 
\end{align}

\subsection{Oscillatory integrals and the slow-rise windowing
  function} \label{sec:win_ieq}

Special considerations must be taken into account in order to solve
the system of equations~\eqref{eq:system1} numerically---mainly in
view of the slow decay of the associated integrands
(equation~\eqref{eq:opers}) for a fixed target point $\nex\in\Gamma$
as $\ney\to\infty$. A direct truncation of the integration domain
(i.e., replacement of the integrals in~\eqref{eq:opers} by
corresponding integrals over the domain $\Gamma^\pm(\nex)\cap\{|\nex
|\leq A\}$) yields slow convergence, on account of edge effects, as
the size $A$ of the truncation domain tends to infinity. Relying on a
certain slow-rise windowing technique that smoothly truncates the
integration domain, the proposed approach addresses this
difficulty---and, in fact, it gives rise to a super-algebraically
convergent algorithm. The resulting windowed integral equations can be
subsequently discretized by means of any integral solver, including,
in particular, the Method of Moments~\cite{Rao1982}, or, indeed, any
Nystr\"om, Galerkin or collocation approach. The particular
implementations presented in this paper are based on the high-order
Nystr\"om method described in~\cite[Sec. 3.5]{Colton2013}.

The proposed methodology is based on use of an infinitely smooth
function $w_A(z)$, defined for $z\in \mathbb{R}$, which satisfies the
following properties: (i)~$w_A(z)$ vanishes for $z$ outside the
interval $[-A,A]$; (ii)~$w_A(z)$ equals 1 for $z \in (-\alpha A, \alpha
A)$ for some $\alpha$ satisfying $0<\alpha <1$; (iii)~All of the
derivatives $w_A^{(p)}(z)$ ($p$ a positive integer) vanish at $z= \pm
A$ and $z = \pm \alpha A$; and (iv)~$w_A(z)$ exhibits a ``slow-rise''
from $0$ to $1$ as $|z|$ goes from $|z| = \pm \alpha A$ to $|z|=
A$---in the sense that each derivative of $w_A$, of any given order,
tends to zero everywhere (and, in particular, in the rise intervals
$\alpha A\leq |z|\leq A$) as $A\to\infty$. The windowing function used
in this paper is given by \ifarxiv
\begin{align}
  w_A(z) = 
  \begin{cases}
    1, & s < 0 \\
    \displaystyle \exp{\bigg( -2 \frac{\exp{(-1/|s|^2)}}{|1-s|^2}  \bigg)}, & \displaystyle 0 \le s \le 1, \hspace{0.8cm} s(z) = \frac{|z|-\alpha A}{ A (1-\alpha)} \\
    0, & s > 0
  \end{cases},
\end{align}
\else
\begin{align}
  w_A(z) = 
  \begin{cases}
    1, & s < 0 \\
    \displaystyle \exp{\bigg( -2 \frac{\exp{(-1/|s|^2)}}{|1-s|^2}  \bigg)}, & \displaystyle 0 \le s \le 1  \\
    0, & s > 0
  \end{cases},
\end{align}
where $ s(z) = \frac{|z|-\alpha A}{ A (1-\alpha)}$, \fi but other
choices could be equally suitable~\cite{Bruno2014}. As shown in that
reference and demonstrated by means of a simple example in
Sec.~\ref{simpl_WI}, such windowing functions can be used to evaluate,
with super-algebraic accuracy, certain improper integrals with slowly
decaying oscillatory integrands---like the integrands in
equation~\eqref{eq:opers}; see Remark~\ref{rem_os}.

\begin{remark} \label{rem_os} 
  In view of the asymptotic expressions for the Hankel function, each
  one of the integral kernels involved in the equation
  system~\eqref{eq:system1} can be expressed in the form
  $h(t)\exp{(it)}$ where $t = k |\nex-\ney| = k |z-z'| \sqrt{1 +
    (\frac{x-x'}{z-z'})^2}$ and where $h(t) \sim t^{-1/2}$
  ($t\to\infty$) is a function each one of whose derivatives is
  bounded for all $t>1$~\cite[Sec. 5.11]{Lebedev1965special} (see
  also~\cite{Demanet:2010cc,Bruno2016}). On the other hand, the
  scattered densities $\varphi^\scat$ approach oscillatory asymptotic
  functions as $\ney = (z',x')$ tends to infinity along
  $\Gamma^\pm(\nex)$~\cite{Nosich1994}: $\varphi^\scat \sim \sum_m
  A^\scat_m u^m_\perp(x') e^{-i k_z^m |z'|}$ as $z' \rightarrow -
  \infty$ and $\varphi^\scat \sim \sum_n B^\scat_n u^n_\perp(x') e^{i
    k_z^n |z'|} $ as $z' \rightarrow \infty$, with similar expressions
  for the density $\psi^\scat$. Combining the kernel and density
  asymptotics it follows, in particular, that the net integrands in
  the operators~\eqref{eq:opers} equals a sum of slowly decaying
  oscillatory functions with wavenumbers $\pm(k + k_z^n)$.
\end{remark}

\subsection{Error estimates for a simplified windowed integral\label{simpl_WI}}
In order to illustrate the properties of the windowed-integration
method used in this paper it is useful to consider here a simple
integration problem presented in~\cite{Bruno2014}, namely, the problem
of numerical evaluation of the integral
\begin{equation}
  \displaystyle I = \int_1^\infty \frac{\e^{i a z}}{\sqrt{z}} dz;
\end{equation}
see, in particular,~\cite[theorem 3.1]{Bruno2014}. Letting
$ I_\text{tr}(A) = \int_1^A \frac{\e^{i a z}}{\sqrt{z}} dz$, $I_w(A) =
\int_1^A w_A(z)\frac{\e^{i a z}}{\sqrt{z}} dz,$ then, by definition $I
= \lim_{A\rightarrow\infty}I_\text{tr}(A)$. As is known, the value of
$I$ is finite---in spite of the slow decay of the integrand (the
integral of the function $1/\sqrt{z}$ in the same domain is
infinite!). The finiteness of the improper integral between $1$ and
$\infty$, which results from cancellation of positive and negative
contributions arising from the oscillatory factor $\e^{i a z}$, may be
verified by integrating by parts the integral $I_\text{tr}(A)$
(differentiating $1/\sqrt{z}$ and integrating $\e^{i a z}$). This procedure
produces two terms: (i)~An integral with a more rapidly decaying integrand and
whose convergence does not rely on cancellations, as well as (ii)~Boundary
contributions at $z=1$ and $z=A$. Besides establishing the existence
of the the limit $\lim _{A\rightarrow\infty}I_\text{tr}(A)$, this
expression tells us that the $z=A$ boundary contribution
$1/(ia\sqrt{A})\e^{i a A}$ equals the error in the approximation of
$I$ by $I_\text{tr}(A)$. On the other hand, use of integration by
parts on $I_w(A)$ does not give rise to a boundary contribution for
$z=A$---on account of the fact $w_A(A)= w_A'(A)=0$. In fact, since all
the derivatives of $w_A(z)$ vanish at $z=A$, the integration by parts
procedure can be performed on $I_w(A)$ an arbitrary number $p$ of
times without ever producing a boundary contribution---a fact which
lies at the heart of the accuracy resulting from the slow-rise
windowing approach. As shown in~\cite[theorem 3.1]{Bruno2014} this
procedure leads to the error estimates \ifarxiv
\begin{equation} 
  \displaystyle |I - I_\text{tr}(A)| = \mathcal{O} \left( \frac{1}{a \sqrt{A}} \right),
\end{equation}
\begin{equation}
  \displaystyle |I - I_w(A)| = \mathcal{O} \left( \frac{1}{\sqrt{a} (a A)^{p-\frac{1}{2}}} \right) \; \text{for every } p \ge 1,
  \label{eq:errest}
\end{equation}
which are corroborated by the results in Table~\ref{table:simple_converge}.
\else
\begin{equation} \label{eq:errest}
  \displaystyle |I - I_\text{tr}(A)| = \mathcal{O} \left( \frac{1}{a \sqrt{A}} \right),~~
  \displaystyle |I - I_w(A)| = \mathcal{O} \left( \frac{1}{\sqrt{a} (a A)^{p-\frac{1}{2}}} \right),
\end{equation}
for every $p \ge 1$. These estimates are corroborated by the results in
Table~\ref{table:simple_converge}.  
\fi
\ifarxiv
\begin{table}[!h]
  \centering
  \begin{tabular}{|c|c|c|}
    \hline
    A & $\displaystyle |I - I_\text{tr}(A)|$ & $\displaystyle |I - I_w(A)|$ \\
    \hline
    10  & $5.0 \times 10^{-2}$ & $3.6 \times 10^{-3}$ \\
    20  & $3.6 \times 10^{-2}$ & $5.8 \times 10^{-5}$ \\
    25  & $3.2 \times 10^{-2}$ & $9.3 \times 10^{-6}$ \\
    50  & $2.2 \times 10^{-2}$ & $3.1 \times 10^{-9}$ \\
    75  & $1.8 \times 10^{-2}$ & $3.1 \times 10^{-12}$ \\
    100 & $1.6 \times 10^{-2}$ & $1.0 \times 10^{-14}$ \\
    \hline
  \end{tabular}
  \caption{Convergence of the windowed integrals $I_\text{tr}(A)$ and $I_w(A)$ with $a=2 \pi$ and window parameter $\alpha=0.5$ \label{table:simple_converge}}
\end{table}
\else
\begin{table}[!h]
  \caption{Convergence of the windowed integrals $I_\text{tr}(A)$ and $I_w(A)$ with $a=2 \pi$ and window 
    parameter $\alpha=0.5$}
  \label{table:simple_converge}
  \centering
  \begin{tabular}{|c|c|c|}
    \hline
    A & $\displaystyle |I - I_\text{tr}(A)|$ & $\displaystyle |I - I_w(A)|$ \\
    \hline
    10  & $5.0 \times 10^{-2}$ & $3.6 \times 10^{-3}$ \\
    20  & $3.6 \times 10^{-2}$ & $5.8 \times 10^{-5}$ \\
    25  & $3.2 \times 10^{-2}$ & $9.3 \times 10^{-6}$ \\
    50  & $2.2 \times 10^{-2}$ & $3.1 \times 10^{-9}$ \\
    75  & $1.8 \times 10^{-2}$ & $3.1 \times 10^{-12}$ \\
    100 & $1.6 \times 10^{-2}$ & $1.0 \times 10^{-14}$ \\
    \hline
  \end{tabular}
\end{table}
\fi

\subsection{Windowed integral equations \label{OW_WGF}}

Per Remark~\ref{rem_os}, for a fixed $\nex \in \Gamma$ the integrands
associated with the operators on the left-hand side of
equation~\eqref{eq:system1} are oscillatory and slowly decaying at
infinity---just like the integrand in the simplified example presented
in Sec.~\ref{simpl_WI}. Thus, it is expected that use of windowing in
the integrands of the aforementioned operators should result in
convergence properties analogous to the ones described in that
section. Since the integrands on the right-hand side (RHS) of
equation~\eqref{eq:system1} are known functions, and since, as shown
below, the full RHS can be evaluated efficiently (by relying on
equation~\eqref{eq:perp}), the use of windowing may be restricted to
the left-hand operators.

These considerations lead to the following system of ``windowed''
integral equations on the {\em bounded} domain $\widetilde \Gamma =
\Gamma \cap \{w_A(\nex) \neq 0 \}$:
\ifarxiv
\begin{equation} \label{eq:system2} 
  E(\nex) \nephi_w^\scat (\nex) + T[w_A\nephi_w^\scat] (\nex)
  = - E(\nex) \nephi^\inc (\nex) - T[\nephi^\inc] (\nex), \quad
  \text{for } \nex \in \widetilde \Gamma.
\end{equation}
\else
\begin{equation} \label{eq:system2} 
  E(\nex) \nephi_w^\scat (\nex) + T[w_A\nephi_w^\scat] (\nex)
  = - E(\nex) \nephi^\inc (\nex) - T[\nephi^\inc] (\nex), 
\end{equation}
for $\nex \in \widetilde \Gamma$.  \fi A discrete version of
equations~\eqref{eq:system2} can be obtained by substituting all
left-hand side integrals by adequate quadrature rules (as mentioned
above, the Nystr\"om method~\cite[Sec. 3.5]{Colton2013} is used for
this purpose in this paper). The windowing function $w_A(\nex)$ used
here is selected as follows: (i)~$w_A=1$ on any portion of $\Gamma$
that is not part of a SIW; (ii)~$w_A$ varies smoothly, with values
between 0 and 1, along any portion of $\Gamma$ contained in a SIW, in
a manner akin to that inherent in the windowing function used in
Sec.~\ref{simpl_WI}, and; (iii) on each SIW, the center of the
$\{w_A(\nex)=1\}$ region lies at the edge of the SIW (see
Fig.~\ref{fig:wgf}). As demonstrated in Sec.~\ref{numer} through a
variety of numerical results, the solution $\nephi_w^\scat$ of
equation~\eqref{eq:system2} provides a super-algebraically accurate
approximation to $\nephi^\scat$ throughout the region $\{w_A(\nex) =
1\}$.

While in principle the integrals $T[\nephi^\inc](\nex)$ on the RHS in
equation~\eqref{eq:system2} could be computed on the basis of windowing
functions, an alternative numerical approach was found preferable. To introduce
this alternative strategy consider the discussion leading to
equation~\eqref{eq:system1}. Clearly, the first (resp. second) component in the
quantity $(\xi,\eta)^t=T[\nephi^\inc](\nex)$ for $\nex \in \widetilde \Gamma$ is
given by the limit as $\nex$ tends to $\widetilde \Gamma$ of the field $u$
(resp. the normal derivative of the field $u$) that results as the pair
$(\phi,\psi)^t$ in~\eqref{eq:rep} is replaced by the incident densities
$(\phi^\inc,\psi^\inc)^t(=\nephi^{\inc})$. Note that, per
equation~\eqref{eq:densities2}, $\nephi^{\inc}$ vanishes identically outside
$\gammainc = \Gamma_j \cap \Omega^\inc$. Now, given that the illuminating
structure is a single SIW whose optical axis coincides with the $z-$axis (as
indicated in Sec.~\ref{sec:ieqf}), the corresponding integration domain
$\gammainc$ can be decomposed as the union
$ \gammainc = \gammainctilde \cup \gammainf$ of the two disjoint segments
$\gammainctilde = \gammainc \cap \{ w_A(\nex) \neq 0 \}$ and
$\gammainf = \gammainc \cap \{ w_A(\nex) = 0 \}$.  (With reference to
Fig.~\ref{fig:wgf} note that
$\Gamma^\inc_j = \left(\bigcup_{\ell=1}^{j-1}\Gamma_{\ell j}^\inc\right) \bigcup
\left(\bigcup_{\ell=j+1}^N\Gamma_{j\ell}^\inc\right) $
where $\Gamma_{j\ell}^\inc = \Gamma_{j\ell} \cap \Omega^\inc$; similarly,
$\gammainctilde$ and $\gammainf$ are decomposed as the unions of the curves
$\widetilde{\Gamma}_{j\ell}^\inc = \Gamma_{j\ell}^\inc \cap \{ w_A(\nex) \neq 0
\}$
and $\Gamma_{j\ell}^\infty = \Gamma_{j\ell}^\inc \cap \{ w_A(\nex) = 0 \}$,
respectively. Fig.~\ref{fig:wgf} only displays the independent components
$\widetilde{\Gamma}_{j\ell}^\inc$ and $\Gamma_{j\ell}^\infty$.)  Using the fact
that the incident densities $\varphi^\inc$ and $\psi^\inc$ vanish outside
$\gammainc$, for $\nex \in \Omega_j$ the relation \ifarxiv
\begin{align} \label{eq:inf1}
  \displaystyle \mathscr{D}[\beta_j \varphi^\inc](\nex) - 
  \mathscr{S}\left[\beta_j \nu_j^{-1}\psi^\inc \right](\nex) = 
  \int_{\gammainctilde \cup \gammainf}  \left[ \frac{\partial G_j}{\partial \nn(\ney)}(\nex, \ney)
  \varphi^\inc(\ney) - G_j(\nex, \ney) \nu_j^{-1}\psi^\inc(\ney) 
  \right] \beta_j \de s_{\ney},
\end{align}
\else
\begin{align} \label{eq:inf1}
  \displaystyle \mathscr{D}&[\beta_j \varphi^\inc] - 
  \mathscr{S}\left[\beta_j \nu_j^{-1}\psi^\inc \right] =  
  \int_{\gammainctilde \cup \gammainf} \beta_j \nonumber \\
  &\displaystyle \times\left[ \frac{\partial G_j}{\partial \nn(\ney)}(\nex, \ney)
  \varphi^\inc(\ney) - G_j(\nex, \ney) \nu_j^{-1}\psi^\inc(\ney) 
  \right]  \de s_{\ney},
\end{align}
\fi results. The evaluation of the right-hand integral is discussed in
what follows.

The integral over the {\em bounded} curve $\gammainctilde$
in~\eqref{eq:inf1} can be computed using any of the standard
singular-integration techniques mentioned in
Sec.~\ref{sec:win_ieq}. The curve $\gammainf$ extends to infinity, on
the other hand, with an integrand that decays slowly (see
Remark~\ref{rem_os}). Fortunately, however, the $\gammainf$
integration problem can be significantly simplified by relying on the
fact that $(\phi^\inc,\psi^\inc)^t(=\nephi^{\inc})$ actually coincides
with the boundary values of the incident field and its normal
derivative, as indicated in equation~\eqref{eq:densities2}. To
evaluate the $\gammainf$ integral consider the identity \ifarxiv
\begin{align}
  \int_{\gammainf} & \left[ \frac{\partial G_j}{\partial \nn(\ney)}(\nex, \ney)
  \varphi^\inc(\ney) - G_j(\nex, \ney) \nu_j^{-1}\psi^\inc(\ney) 
  \right] \beta_j \de s_{\ney}  = \nonumber \\
  &-\int_{\Gamma_j^\perp}  \left[ \frac{\partial G_j}{\partial \nn(\ney)}(\nex, \ney)
  u_+^\inc(\ney) - G_j(\nex, \ney) \nu_j^{-1}\frac{\p u_+^\inc}{\p \nn}(\ney) 
  \right] \beta_j \de s_{\ney},  \quad \text{for }\nex \in \{z>-A\},
   \label{eq:perp}
\end{align}
\else
\begin{align}
  \int_{\gammainf}  \beta_j \bigg[& \left.\frac{\partial G_j}{\partial \nn(\ney)}(\nex, \ney)
  \varphi^\inc(\ney) - G_j(\nex, \ney) \nu_j^{-1}\psi^\inc(\ney) 
  \right] \de s_{\ney}   \nonumber \\
  =&-\int_{\Gamma_j^\perp}  \beta_j \left[ \frac{\partial G_j}{\partial \nn(\ney)}(\nex, \ney)
  u_+^\inc(\ney) \right. \nonumber \\ & \left. - G_j(\nex, \ney) \nu_j^{-1}\frac{\p u_+^\inc}{\p \nn}(\ney) 
  \right] \de s_{\ney} , \quad \nex \in \{z>-A\}
   \label{eq:perp}
\end{align}
\fi that results by using Green's theorem in a manner akin to~\cite{DeSanto1997}
(the corresponding bounded-domain result can be found e.g. in~\cite[theorem
3.1]{Colton1983}).  (Here $\Gamma_j^\perp = \Omega_j \cap \{ z=-A\}$ is such
that $\Gamma_j^\perp \cup \gammainf$ is the boundary of the region
$\Omega_j \cap \{z<-A\}$, see Fig.~\ref{fig:wgf}.) Equation~\eqref{eq:perp}
provides an alternative approach for the evaluation of the integrals over
$\gammainf$ in equation~\eqref{eq:inf1}. While some of the $\Gamma_j^\perp$
curves are unbounded, along such unbounded curves the fields decay
\emph{exponentially fast} (see equation~\eqref{eq:mode}), and, thus, arbitrary
accuracy can efficiently be achieved in the corresponding integrals by
truncation of the integration domain to a relatively small bounded portion of
$\Gamma_j^\perp$.

\ifarxiv
\begin{figure}[h!] \centering 
  \subfloat[Loglog plot] {\includegraphics[width=3.3in]{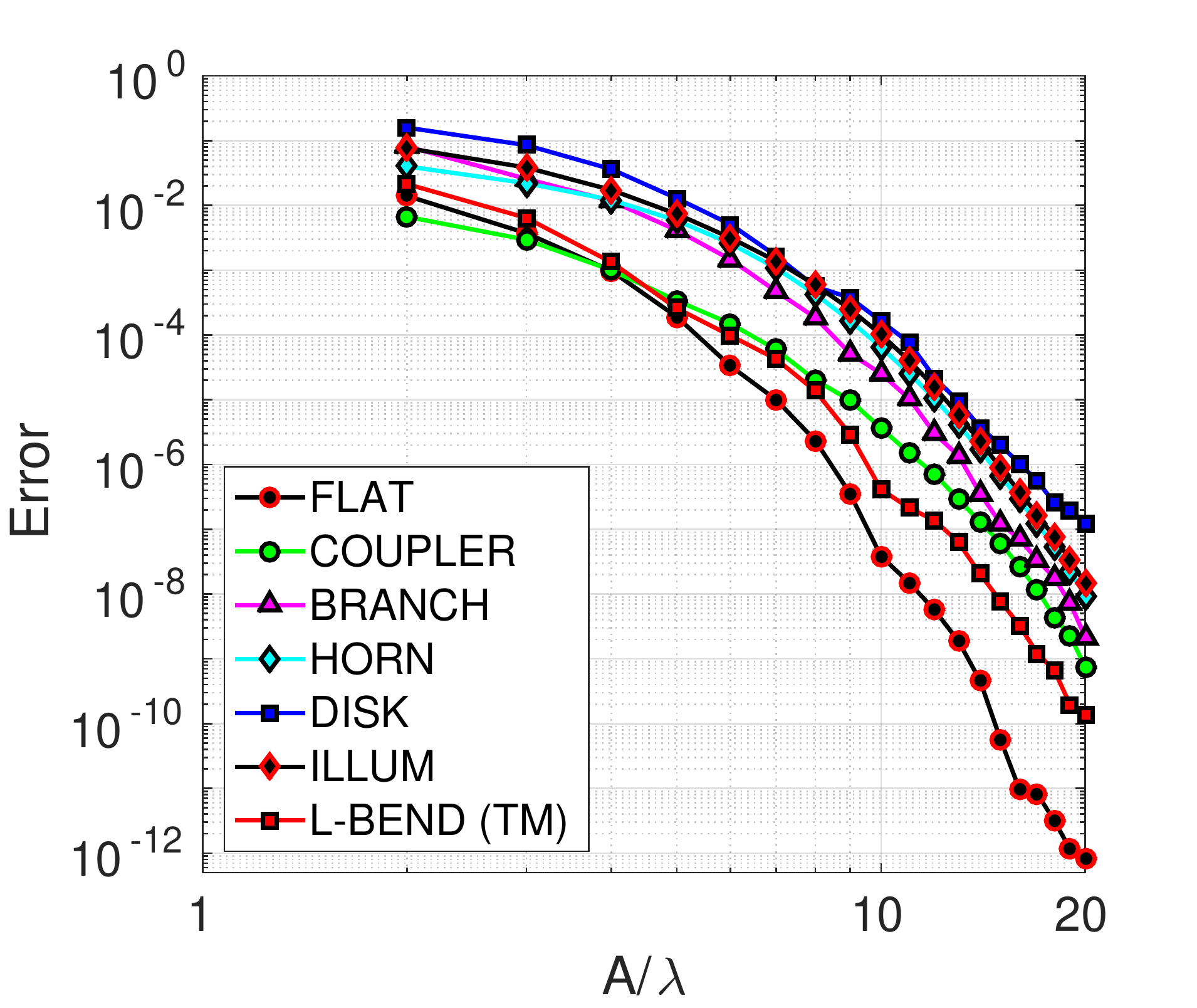} 
  \label{fig:log}}  
  \subfloat[Semilog plot] {\includegraphics[width=3.3in]{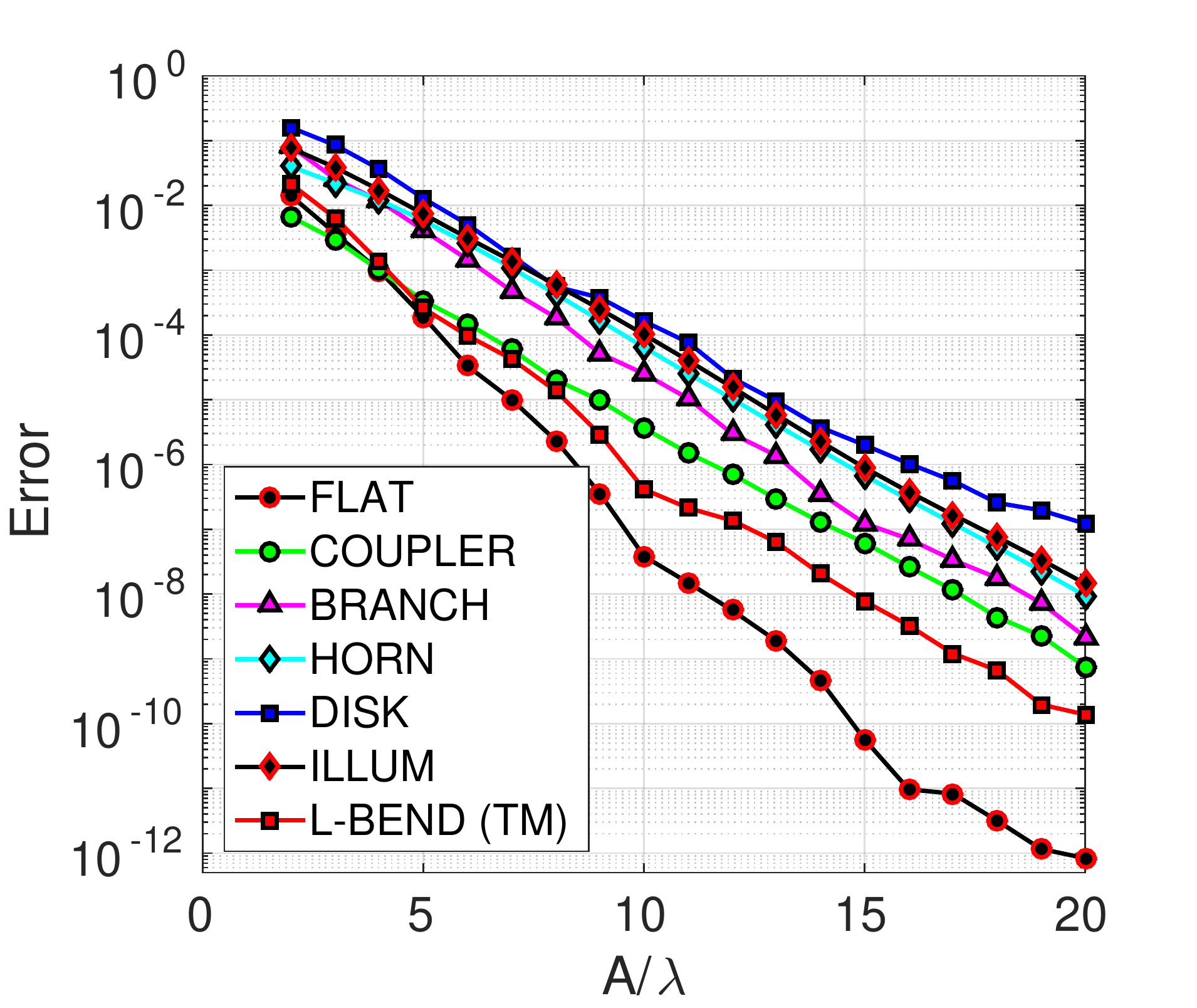}%
    \label{fig:semi}} 
   \caption{Super-algebraic convergence of the WGF method for various
    test configurations. In all cases a sufficiently fine numerical
    grid was used (eight to twelve points per wavelength) to ensure
    the leading error is caused by the WGF's slow-rise boundary
    truncation. Here the scaled window size $A/\lambda$ is varied
    while $\lambda = \text{max}_j(2 \pi / k_j)$ is kept fixed (where
    the values of $k_j$ are described in the text for each test
    case).}
  \label{fig:convergence}
\end{figure}
\else
\fi

The overall WGF method for open waveguides is summarized in points~1
to~3 below. As mentioned above in this paper, the bounded-interval
numerical integrations mentioned in the points~1 to~3 can be effected
by means of any numerical integration method applicable to the kinds
of singular integrals (with singular kernels and, at corners, unknown
singular densities) associated with the problems under
consideration. The implementation used to produce the numerical
results in Sec.~\ref{numer} relies on numerical integration methods
derived in a direct manner from those described
in~\cite[Sec. 3.5]{Colton2013}.

The WGF algorithm proceeds as follows:
\begin{enumerate}
\item[(1)] Evaluate the RHS of equation~\eqref{eq:system2} by decomposing the
  integrals in the operators~\eqref{eq:oper} into integrals over
  $\gammainctilde$ and $\gammainf$. The integrals involving the bounded
  integration domains $\gammainctilde$ are computed by direct numerical
  integration. Relying on the exponential decay of the integrands of the RHS of
  equation~\eqref{eq:perp}, on the other hand, the integrals over $\gammainf$
  are obtained by truncation of the integrals over $\Gamma_j^\perp$.
\item[(2)] Solve for $\nephi^\scat_w$ in equation~\eqref{eq:system2} by
  either inverting directly the corresponding discretized system (as
  is done in this paper), or by means of a suitable iterative
  linear-algebra solver, if preferred.
\item[(3)] Evaluate the approximate fields $u_w$ using the
  representation formula~\eqref{eq:rep} in conjunction with
  equation~\eqref{eq:perp}:
\ifarxiv
\begin{align}
    u_w(\nex) =& \mathscr{D}[\beta_j w_A \varphi_w^\scat](\nex) - 
  \mathscr{S}\left[\beta_j \nu_j^{-1}w_A\psi_w^\scat \right](\nex)  + \mathscr{D}[\beta_j \varphi^\inc](\nex) - 
  \mathscr{S}\left[\beta_j \nu_j^{-1}\psi^\inc \right](\nex),
  \end{align}
\else
  \begin{align}
    u_w(\nex) =& \mathscr{D}[\beta_j w_A \varphi_w^\scat](\nex) - 
  \mathscr{S}\left[\beta_j \nu_j^{-1}w_A\psi_w^\scat \right](\nex) \nonumber \\
  +& \mathscr{D}[\beta_j \varphi^\inc](\nex) - 
  \mathscr{S}\left[\beta_j \nu_j^{-1}\psi^\inc \right](\nex),
  \end{align} \fi
  for $\nex \in \Omega_j$, and
  where the layer potentials involving $\varphi^\inc$ and $\psi^\inc$
  are computed as in point~1, that is, by direct evaluation of the
  integrals over $\gammainctilde$ and via equation~\eqref{eq:perp} for
  the integrals over $\gammainf$.
\end{enumerate}
It can be seen~\cite{Bruno2016} that, as $A\to\infty$, the fields in
point~3 are super-algebraically accurate in any bounded region in the
plane.

\ifarxiv
\begin{figure}[h!] \centering  
  \subfloat[COUPLER] {\includegraphics[width=2.5in]{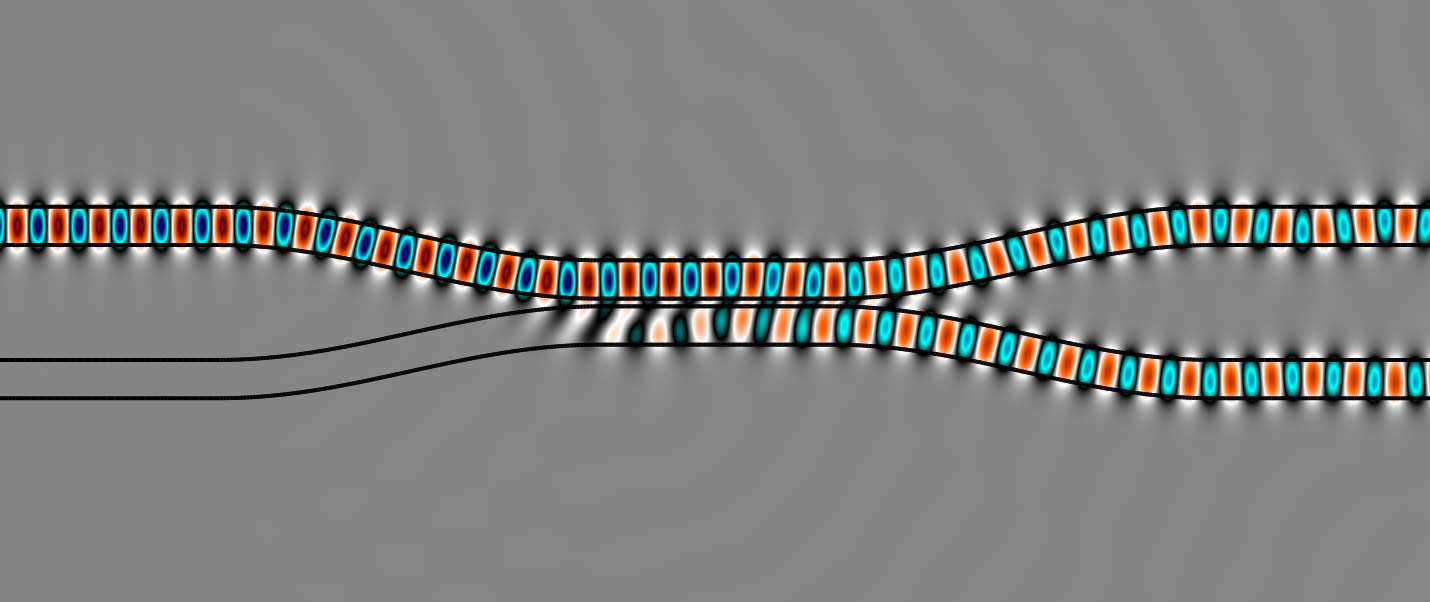} 
  \label{fig:coupler_re}}  
  \subfloat[COUPLER] {\includegraphics[width=2.5in]{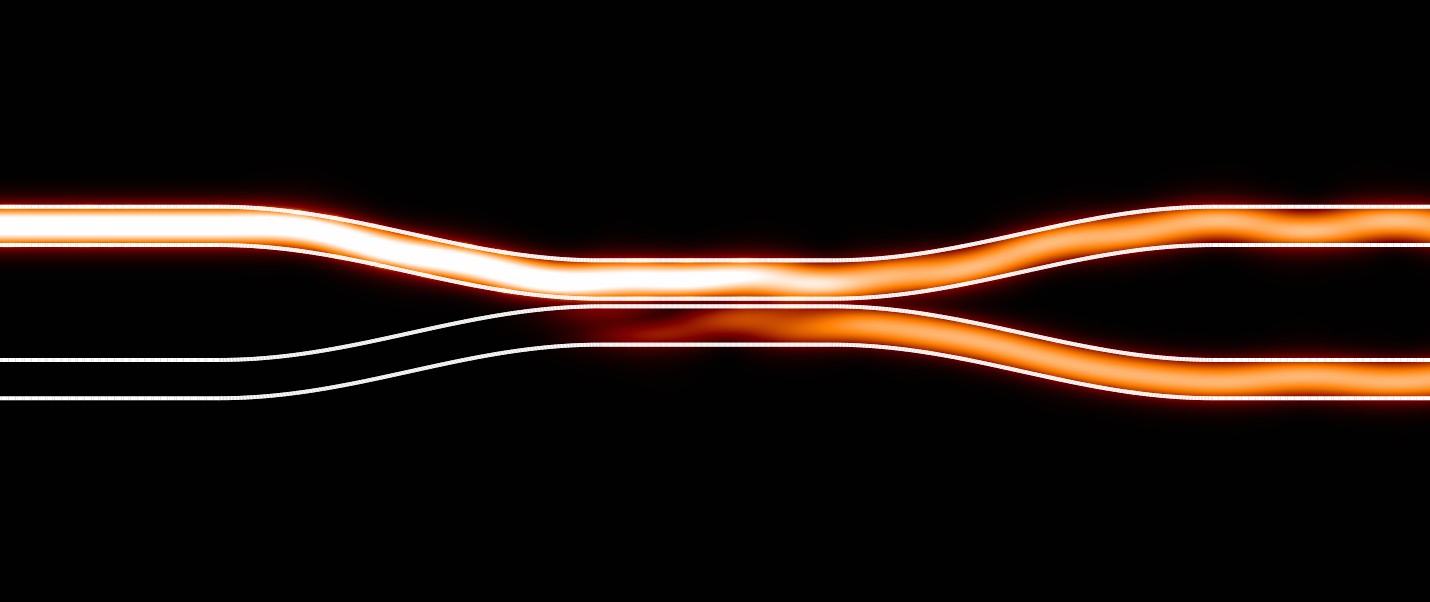}%
    \label{fig:coupler_abs}} \\
  \subfloat[BRANCH] {\includegraphics[width=2.5in]{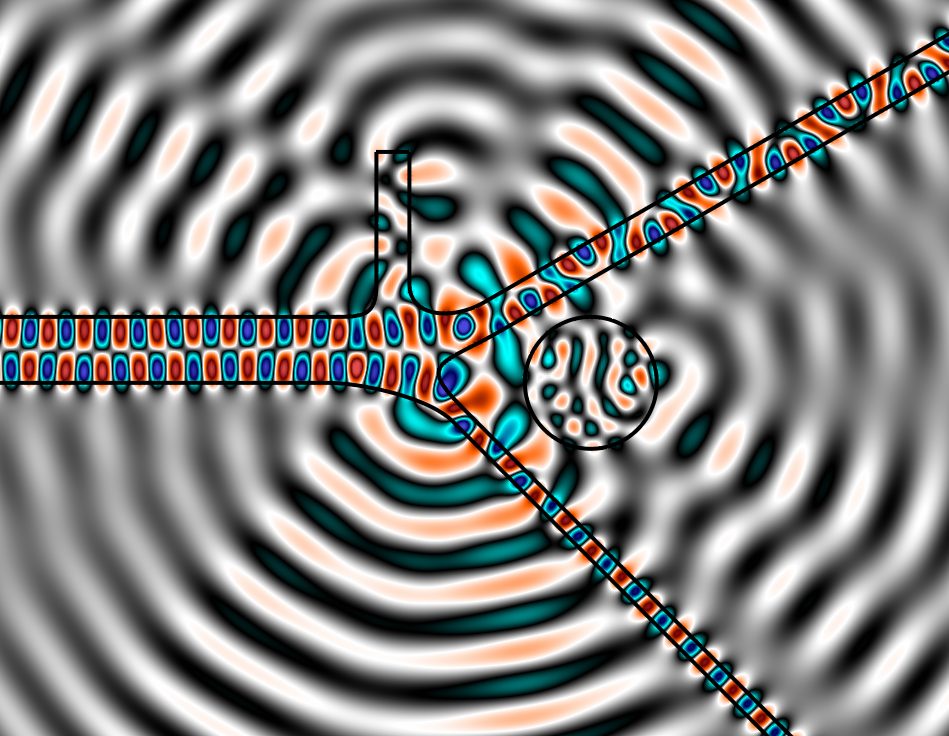} 
  \label{fig:branches_re}} 
  \subfloat[BRANCH] {\includegraphics[width=2.5in]{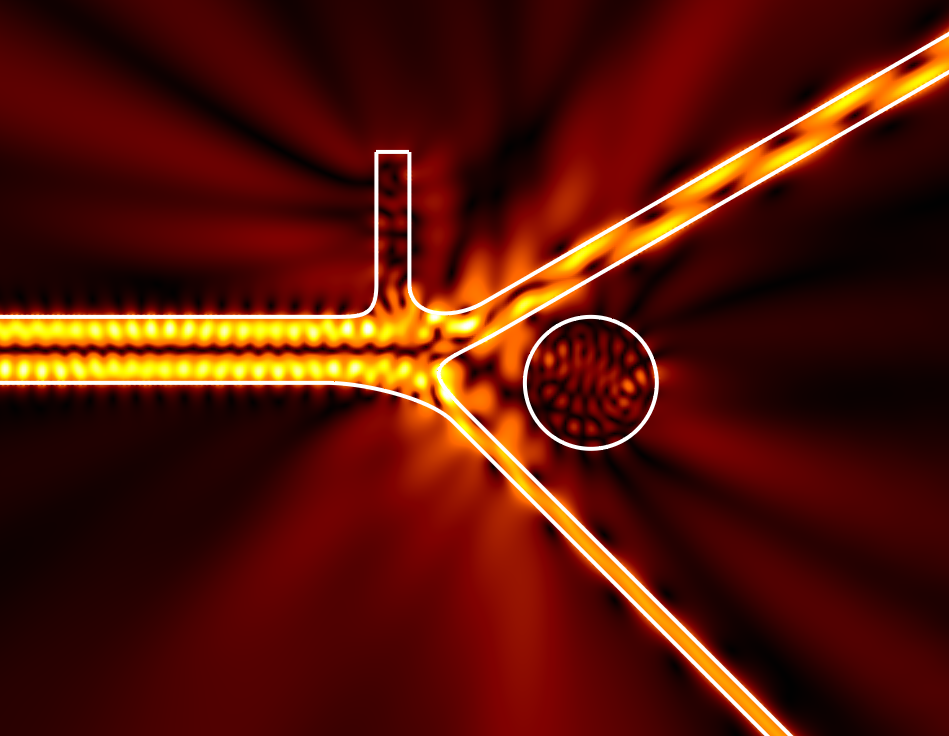}%
    \label{fig:branches_abs}} \\ 
  \subfloat[HORN] {\includegraphics[width=2.5in]{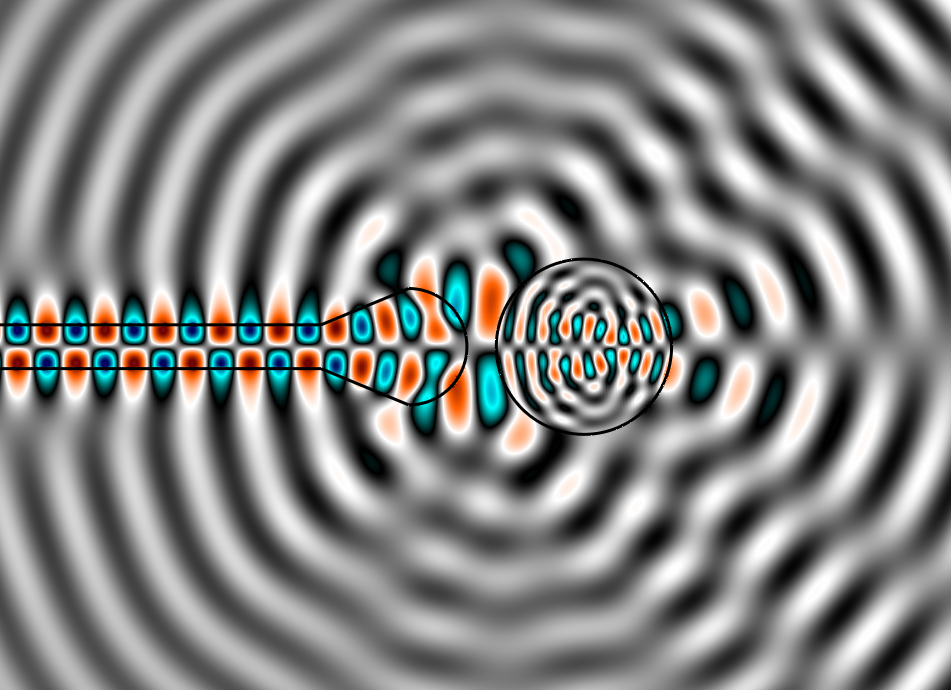} 
  \label{fig:horn_re}} 
  \subfloat[HORN] {\includegraphics[width=2.5in]{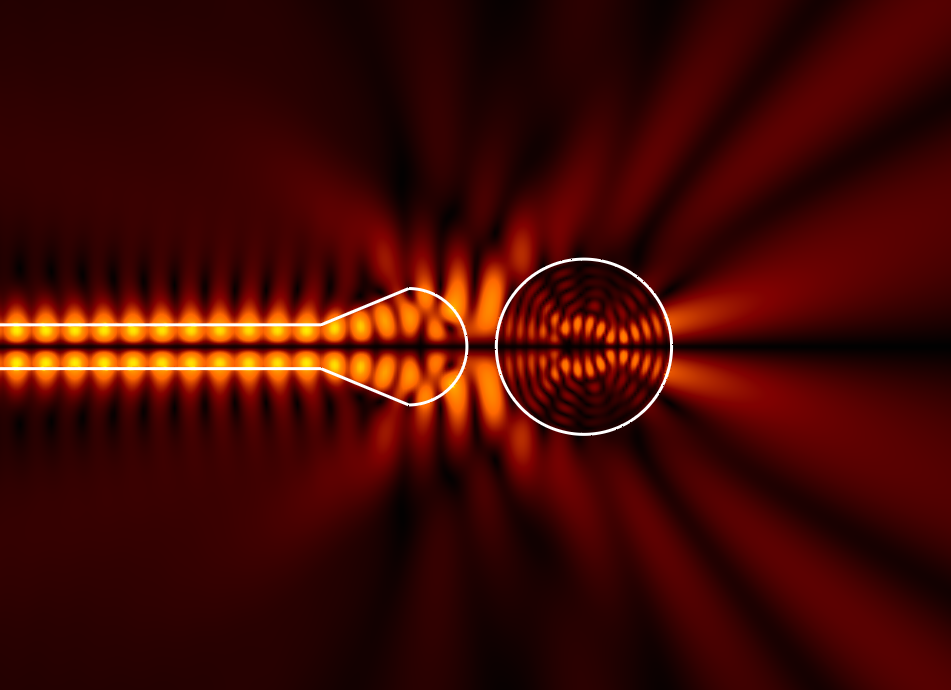}%
     \label{fig:horn_abs}}  \\ 
  \subfloat[Colormap for $\text{Re}(u)$] {\includegraphics[width=2.5in]{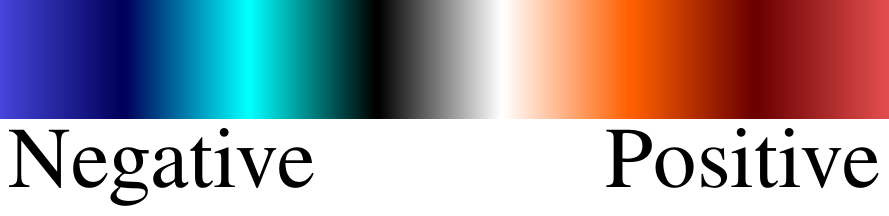} 
    \label{fig:colormap_re}} 
  \subfloat[Colormap for $|u|$] {\includegraphics[width=2.5in]{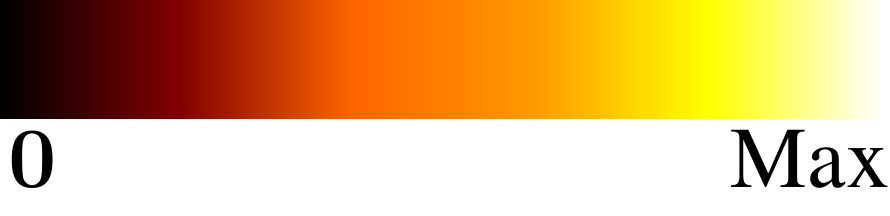}%
    \label{fig:colormap_abs}} 
  \caption{Real part and absolute value of $u_w$ (left and right
    columns, respectively) produced by the WGF method for several
    open-waveguide problems.  }
  \label{fig:fields1}
\end{figure}
\else
\fi

\section{Numerical examples\label{numer}}

This section demonstrates the character of the WGF method introduced
in Sec.~\ref{sec:wgf} through a variety of numerical results. These
results were obtained by means of a Matlab implementation of the
algorithms described in Sec.~\ref{OW_WGF} on a six-core 3.40 GHz Intel
i7-4930K Processor with 12 Mb of cache and 32 Gb of RAM. As mentioned
in Sec.~\ref{sec:win_ieq}, a Nystr\"om method was utilized, where the
number of points per wavelength was selected in such a way that the
dominant error arises from the windowed truncation. The results in
this paper were obtained by means of discretizations containing a
number of eight to twelve points per wavelength, as needed to guarantee the
accuracy reported for each numerical solution. The reported errors
were evaluated by means of the expression
\begin{align} \label{eq:error}
  \displaystyle \text{Error} = \sqrt{ \frac{ \sum_{i=1}^{N_d} |u_w(\nex_i)-u^\text{ref}(\nex_i)|^2}{
  \sum_{i=1}^{N_d} |u^\text{ref}(\nex_i)|^2}}
\end{align}
with $N_d=100$, and where the $N_d$ points $\nex_i$ lie along certain
curves, selected in different manner for each test case, along which
significant features of the numerically computed fields were observed.
The errors in all the nonuniform waveguide problems were evaluated by
means of equation~\eqref{eq:error} with reference solutions
$u^\text{ref}$ produced by means of window functions $w_A$ with
$A/\lambda > 20$.

The first numerical example considered here is a simple uniform
waveguide---for which, physically, a single mode
(equation~\eqref{eq:mode}), or even a superposition of such modes, can
propagate from $-\infty$ to $+\infty$ without disturbance. This simple
problem provides a direct test for the accuracy of the proposed WGF
approach, for which errors were evaluated by using exact solutions as
references: $u^\text{ref}=u$ in this case. The incident field was
prescribed on the region $\Omega^\inc=\{z\le0\}$.  The observed error
as a function of the window size is displayed in
Fig.~\ref{fig:convergence} under the label ``FLAT''.

\ifarxiv
\begin{figure}[h!] \centering  
  \subfloat[DISK] {\includegraphics[width=2.5in]{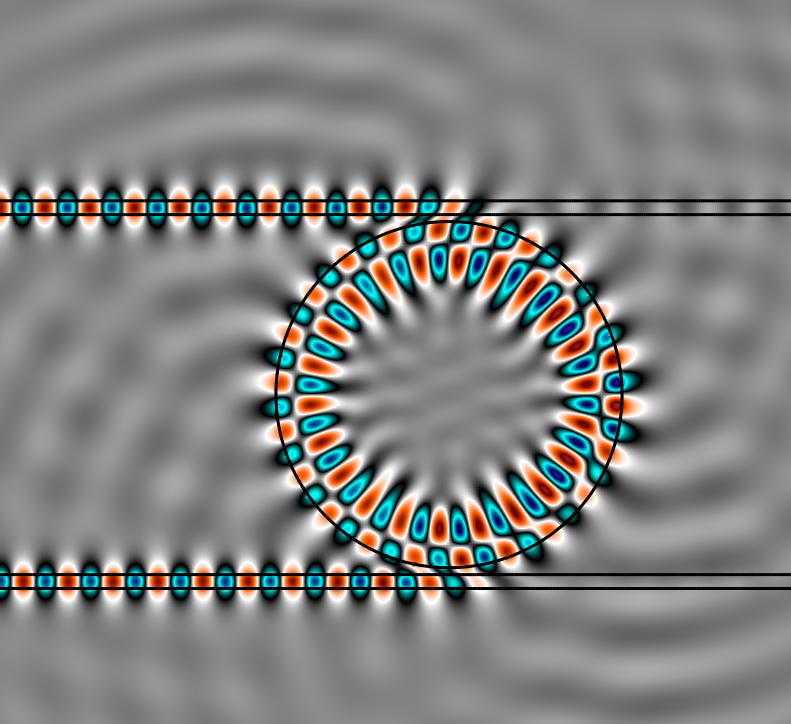} 
    \label{fig:disk_re}} 
  \subfloat[DISK] {\includegraphics[width=2.5in]{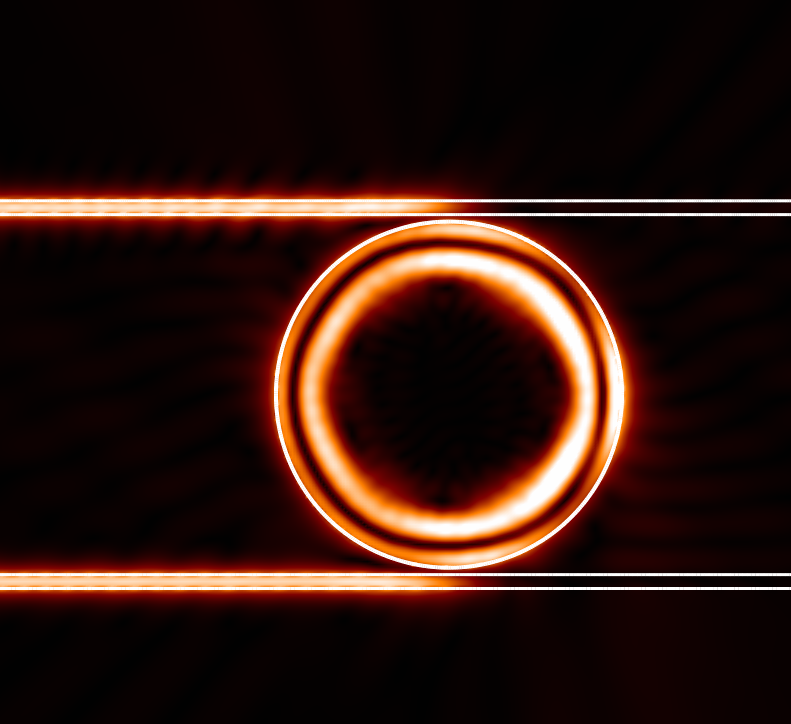}%
     \label{fig:disk_abs}}  \\
  \subfloat[ILLUM] {\includegraphics[width=2.5in]{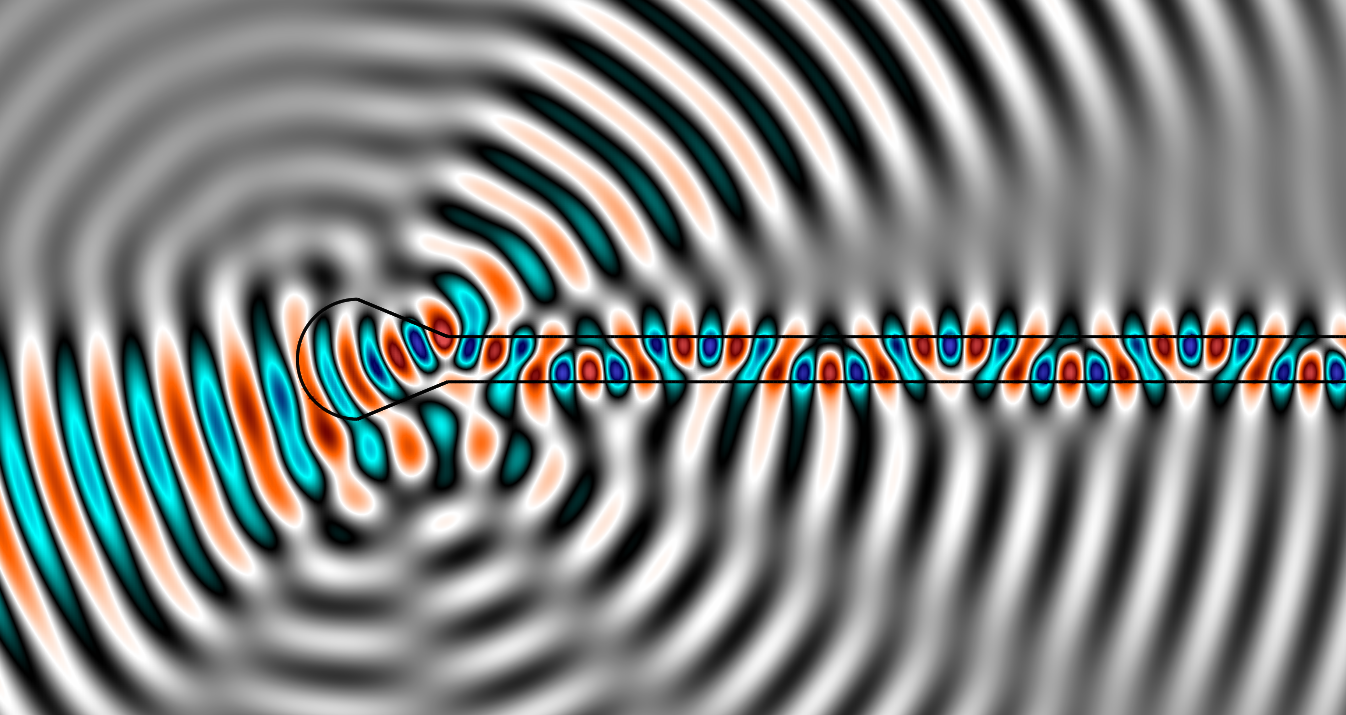} 
    \label{fig:illumination_re}} 
  \subfloat[ILLUM] {\includegraphics[width=2.5in]{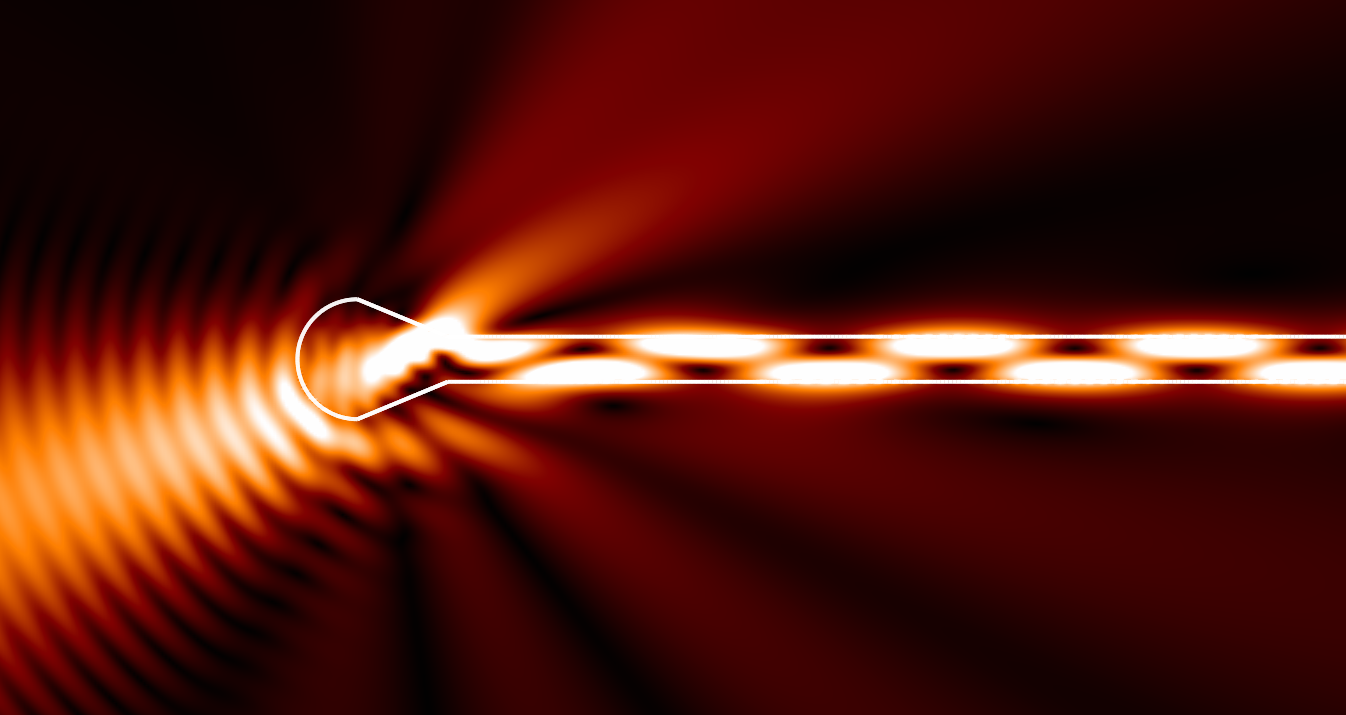}%
    \label{fig:illumination_abs}} \\
  \subfloat[L-BEND] {\includegraphics[width=2.5in]{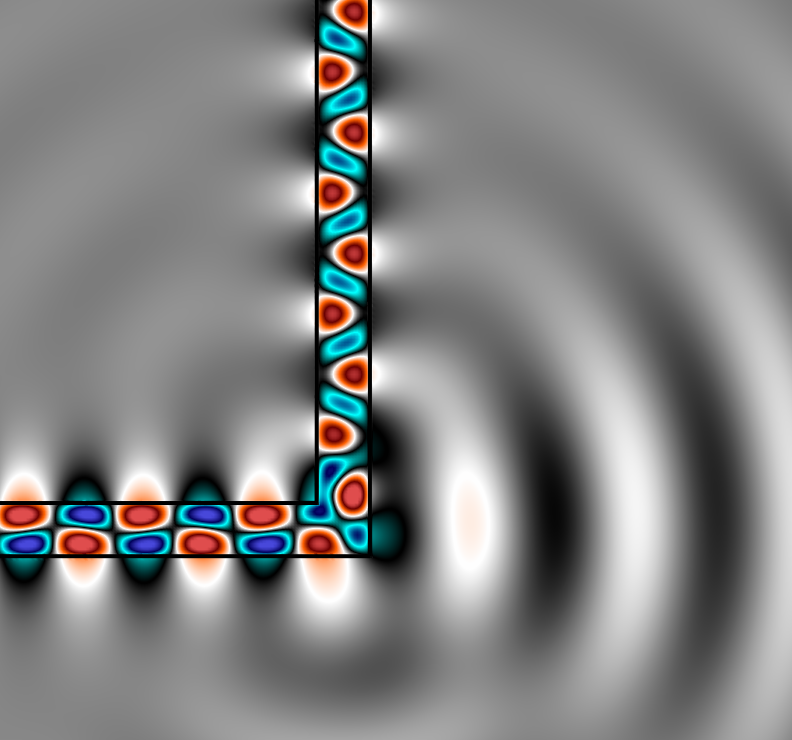} 
    \label{fig:l-shaped_re}} 
  \subfloat[L-BEND] {\includegraphics[width=2.5in]{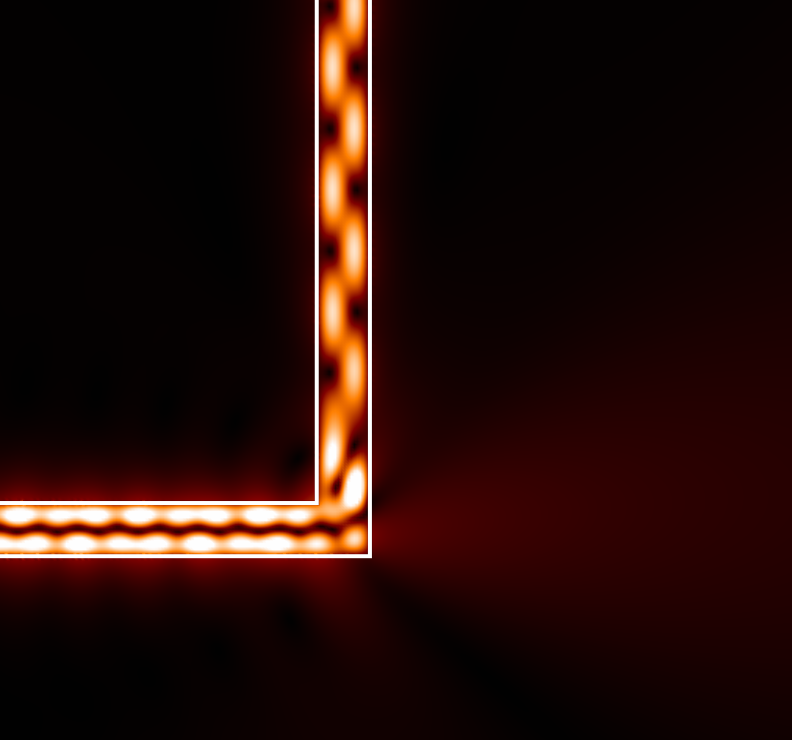}%
    \label{fig:l-shaped_abs}} \\
  \subfloat[Colormap for $\text{Re}(u)$] {\includegraphics[width=2.5in]{colormap_real.pdf} 
    \label{fig:colormap_re}} 
  \subfloat[Colormap for $|u|$] {\includegraphics[width=2.5in]{colormap_abs.pdf}%
    \label{fig:colormap_abs}} 
  \caption{Real part and absolute value of $u_w$ (left and right
    columns, respectively) produced by the WGF method for several
    open-waveguide problems.  }
  \label{fig:fields2}
\end{figure}

\else
\fi

Error curves for a number of additional configurations are also
included in Fig.~\ref{fig:convergence}, and corresponding near field
images are displayed in \ifarxiv Figs.~\ref{fig:fields1}
and~\ref{fig:fields2} \else Fig.~\ref{fig:fields}\fi. (The ``trivial''
near field image for the uniform waveguide mentioned above is not
included in \ifarxiv Figs.~\ref{fig:fields1} and~\ref{fig:fields2}
\else Fig.~\ref{fig:fields}\fi.)  The configurations considered (under
TE or TM polarizations, as indicated in each case) are as follows:
\begin{itemize}
\item COUPLER (TE). Optical coupler illuminated by the first symmetric
  mode incoming from the top-left SIW. Waveguide core wavenumber
  $k_\text{co} = 2 \pi$, cladding wavenumber $k_\text{cl}=\pi$ and
  waveguide half-width $h=0.5$ were used. The waveguide separation in
  the cross-talk region was set to $0.2$.
\item BRANCH (TE). Branching waveguide structure illuminated by the
  first antisymmetric mode incoming from the left SIW. The half-widths
  of the horizontal, top-right and bottom right SIWs are $1$, $0.5$
  and $0.25$, respectively, and the angle between the two branching
  waveguides is $5\pi/12$ radians. The bottom-right SIW is a
  single-mode waveguide: only one specific mode can propagate along
  this structure. The wavenumbers in the core, cladding and circular
  obstacle regions are $k_\text{co} = 2\pi$, $k_\text{cl} = \pi$ and
  $k_\text{ob} = 5\pi/2$, respectively.
\item HORN (TE). Terminated waveguide horn antenna illuminating a dielectric circular
  obstacle. The system is illuminated by the first antisymmetric mode of the
  waveguide.  The core and cladding wavenumbers are $k_\text{co}=2\pi$ and
  $k_\text{cl}=4\pi/3$ respectively, and the circular obstacle wavenumber is
  $k_\text{ob}=4\pi$. The half-width of the waveguide is $h=0.5$ and the radius
  of the obstacle is $2$.
\item DISK (TE). Circular disk
  resonator~\cite[Sec. 16.5]{taflove}. The device is illuminated by
  the first symmetric mode incoming from the top-left waveguide. Both
  waveguides have core wavenumbers $k_\text{co}=250 \pi/127$, the
  cladding wavenumber is $k_\text{cl} = 125\pi/127$ and the disk has
  the same wavenumber as the core regions. The waveguide half-widths
  are both $h=0.2$, the disk has radius $5$, and separation between
  the waveguides and the disk is $0.2$. The wavenumbers were selected
  to excite a near resonance in the circular cavity.
\item ILLUM (TE). Excitation of waveguide modes in a terminated waveguide. The
  structure is illuminated by a beam incoming from the left at an angle $\pi/10$
  radians below the horizontal. (A description of the illuminating field is
  presented below.) The core and cladding wavenumbers are given by
  $k_\text{co}=2\pi$ and $k_\text{cl} = 4\pi/3$ respectively, while the
  half-width of the waveguide is $h=0.5$.
\item L-BEND (TM). Sharp $L-$bend illuminated by the first antisymmetric mode
  incoming from the left waveguide. The core and cladding wavenumbers are
  $k_\text{co}=2\pi$ and $k_\text{cl}=2\pi/3$ respectively, and the waveguide
  half-width is given by $h=0.5$.
\end{itemize}
The labels used here (COUPLER, BRANCH, etc.)  correspond to those used
in Figs.~\ref{fig:convergence}\ifarxiv,~\ref{fig:fields1},
and~\ref{fig:fields2}. \else and~\ref{fig:fields}.\fi The computing
times required to achieve an accuracy better than $1.0\times10^{-6}$
are presented in Table~\ref{table:time}.

The illuminating field used for the configuration ILLUM is given by
the angular spectrum representation~\cite[Eq. (49)]{DeSanto1997}
\begin{align}
  u^\inc(r,\theta) = \int_{-\pi/2}^{\pi/2}F(\alpha) e^{i k_1 r \cos{(\theta + \alpha)}} 
  \de \alpha,
\end{align}
with $F(\alpha)=e^{-12.5(\alpha+\pi/10)^2}$. Note that
equation~\eqref{eq:system1} does not directly apply in this case since the
illuminating field used here is not a waveguide mode. But only slight
modifications are necessary: the relevant equation in this case is
\begin{equation} \label{eq:system3} E(\nex) \nephi^\scat (\nex) + T[\nephi^\scat] (\nex)
  = \nephi^\inc (\nex) , \quad \nex \in \Gamma.
\end{equation}
This equation can be solved by a procedure analogous to that presented
in Sec.~\ref{OW_WGF}.

\ifarxiv
\begin{table}[!h]
  \centering
  \begin{tabular}{|l|c|c|c|c|}
    \hline
    Problem & \#-Unknowns & $A/\lambda$ & Time-$\nephi^\scat$ & Time-$u_w$ \\
    \hline
    FLAT    & 1752 & 9  & 1.06  & 3.47 \\
    COUPLER & 6450 & 12 & 13.05 & 9.23 \\
    BRANCH  & 5978 & 14 & 13.82 & 8.13 \\
    HORN    & 2902 & 15 & 3.28  & 7.69 \\
    DISK    & 6556 & 16 & 11.62 & 5.09 \\
    ILLLUM  & 1374 & 15 & 0.83  & 3.34 \\
    L-BEND & 2516 & 12 & 2.03  & 3.80 \\
    \hline
  \end{tabular}
  \caption{
    Computing times (in seconds) required by the WGF method to produce the
    densities $\nephi^\scat$ (on the discretization grid) and to evaluate
    the fields $u_w$ (on a $128\times 128$ evaluation grid), with an
    accuracy better than  $1.0\times10^{-6}$, for the various test cases mentioned in the
    text. \label{table:time}}
\end{table}

\else
\begin{table}[!h]
  \caption{
    Computing times (in seconds) required by the WGF method to produce the
    densities $\nephi^\scat$ (on the discretization grid) and to evaluate
    the fields $u_w$ (on a $128\times 128$ evaluation grid), with an
    accuracy better than  $1.0\times10^{-6}$, for the various test cases mentioned in the
    text.}
  \label{table:time}
  \centering
  \begin{tabular}{|l|c|c|c|c|}
    \hline
    Problem & \#-Unknowns & $A/\lambda$ & Time-$\nephi^\scat$ & Time-$u_w$ \\
    \hline
    FLAT    & 1752 & 9  & 1.06  & 3.47 \\
    COUPLER & 6450 & 12 & 13.05 & 9.23 \\
    BRANCH  & 5978 & 14 & 13.82 & 8.13 \\
    HORN    & 2902 & 15 & 3.28  & 7.69 \\
    DISK    & 6556 & 16 & 11.62 & 5.09 \\
    ILLLUM  & 1374 & 15 & 0.83  & 3.34 \\
    L-BEND & 2516 & 12 & 2.03  & 3.80 \\
    \hline
  \end{tabular}
\end{table}
\fi

\section{Conclusions} \label{sec:conclusions} The WGF method
introduced in this paper provides super-algebraically accurate
approximations as the window sizes $A$ are increased. The proposed
approach retains the attractive qualities of boundary integral
equation methods, such as reduced dimensionality, efficient
parallelization and high-order accuracy for arbitrary geometries. And,
while the present implementation is based on use of Nystr\"om
integral-equation solvers, any available boundary integral method for
transmission problems, such as, e.g., those based on the Method of
Moments, can be easily modified to incorporate the WGF methodology.
The ideas presented here can be generalized to the three-dimensional
case; such generalizations will be presented elsewhere.

\section*{Acknowledgment}
The authors gratefully acknowledge support by NSF and AFOSR through
contracts DMS-1411876 and FA9550-15-1-0043, and by the NSSEFF Vannevar
Bush Fellowship under contract number N00014-16-1-2808.

\ifarxiv
  \bibliography{bibliography}

\begin{thebibliography}{10}

\bibitem{Berenger1994}
J.-P. Berenger.
\newblock {A perfectly matched layer for the absorption of electromagnetic
  waves}.
\newblock {\em Journal of Computational Physics}, 114(2):185--200, Oct. 1994.

\bibitem{Bruno2014}
O.~P. Bruno and B.~Delourme.
\newblock {Rapidly convergent two-dimensional quasi-periodic Green function
  throughout the spectrum-including Wood anomalies}.
\newblock {\em Journal of Computational Physics}, 262:262--290, 2014.

\bibitem{Bruno2016}
O.~P. Bruno, M.~Lyon, C.~P\'{e}rez-Arancibia, and C.~Turc.
\newblock {Windowed Green Function Method for Layered-Media Scattering}.
\newblock {\em SIAM Journal on Applied Mathematics}, 76(5):1871--1898, Sept.
  2016.

\bibitem{royal}
O.~P. Bruno, S.~P. Shipman, C.~Turc, and S.~Venakides.
\newblock Superalgebraically convergent smoothly windowed lattice sums for
  doubly periodic green functions in three-dimensional space.
\newblock {\em Proceedings of the Royal Society of London A: Mathematical,
  Physical and Engineering Sciences}, 472(2191), 2016.

\bibitem{chaubell2009}
J.~Chaubell, O.~P. Bruno, and C.~O. Ao.
\newblock Evaluation of em-wave propagation in fully three-dimensional
  atmospheric refractive index distributions.
\newblock {\em Radio Science}, 44(1), 2009.

\bibitem{Ciraolo2005}
G.~Ciraolo.
\newblock {\em {Non-Rectilinear Waveguides: Analytical and Numerical Results
  Based on the Green's Function}}.
\newblock PhD thesis, Universita Degli Studi Di Firenze, 2005.

\bibitem{Ciraolo2008}
G.~Ciraolo and R.~Magnanini.
\newblock {Analytical results for 2-D non-rectilinear waveguides based on the
  Green's function}.
\newblock {\em Mathematical Methods in the Applied Sciences}, 31(13):1--18,
  2008.

\bibitem{Colton1983}
D.~Colton and R.~Kress.
\newblock {\em {Integral Equation Methods in Scattering Theory}}.
\newblock 1983.

\bibitem{Colton2013}
D.~Colton and R.~Kress.
\newblock {\em {Inverse Acoustic and Electromagnetic Scattering Theory}}.
\newblock Springer, New York, third edition, 2013.

\bibitem{Demanet:2010cc}
L.~Demanet and L.~Ying.
\newblock {Scattering in flatland: efficient representations via wave atoms}.
\newblock {\em Foundations of Computational Mathematics. The Journal of the
  Society for the Foundations of Computational Mathematics}, 10(5):569--613,
  2010.

\bibitem{DeSanto1997}
J.~DeSanto and P.~A. Martin.
\newblock {On the derivation of boundary integral equations for scattering by
  an infinite one-dimensional rough surface}.
\newblock {\em The Journal of the Acoustical Society of America}, 102(1):67,
  1997.

\bibitem{Farjadpour2006}
A.~Farjadpour, D.~Roundy, A.~Rodriguez, M.~Ibanescu, P.~Bermel, J.~D.
  Joannopoulos, S.~G. Johnson, and G.~W. Burr.
\newblock {Improving accuracy by subpixel smoothing in the finite-difference
  time domain}.
\newblock {\em Optics Letters}, 31(20):2972, Oct. 2006.

\bibitem{Lebedev1965special}
N.~N. Lebedev.
\newblock {\em Special functions and their applications}.
\newblock Prentice-Hall, Inc., 1965.

\bibitem{Magnanini2001}
R.~Magnanini and F.~Santosa.
\newblock {Wave Propagation in a 2-D Optical Waveguide}.
\newblock {\em SIAM Journal on Applied Mathematics}, 61(4):1237--1252, 2001.

\bibitem{Monro2007}
J.~A. Monro.
\newblock {\em {A Super-Algebraically Convergent, Windowing-Based Approach to
  the Evaluation of Scattering from Periodic Rough Surfaces}}.
\newblock PhD thesis, California Institute of Technology, 2007.

\bibitem{Nosich1994}
A.~I. Nosich.
\newblock {Radiation conditions, limiting absorption principle, and general
  relations in open waveguide scattering}.
\newblock {\em Journal of electromagnetic waves and applications},
  8(3):329--353, 1994.

\bibitem{Perez-Arancibia2017}
C.~Perez-Arancibia.
\newblock {\em {Windowed integral equation methods for problems of scattering
  by defects and obstacles in layered media}}.
\newblock PhD thesis, California Institute of Technology, 2016.

\bibitem{Rao1982}
S.~Rao, D.~Wilton, and A.~Glisson.
\newblock {Electromagnetic scattering by surfaces of arbitrary shape}.
\newblock {\em IEEE Transactions on Antennas and Propagation}, 30(3):409--418,
  May 1982.

\bibitem{taflove}
A.~Taflove and S.~C. Hagness.
\newblock {\em Computational Electrodynamics: The Finite-Difference Time-Domain
  Method}.
\newblock Artech House, Inc., third edition, 2005.

\bibitem{Zhang2011}
L.~Zhang, J.~H. Lee, A.~Oskooi, A.~Hochman, J.~K. White, and S.~G. Johnson.
\newblock {A Novel Boundary Element Method Using Surface Conductive Absorbers
  for Full-Wave Analysis of 3-D Nanophotonics}.
\newblock {\em Journal of Lightwave Technology}, 29(7):949--959, Apr. 2011.

\end{thebibliography}
\else
\fi

\end{document}